\documentclass{dcdssissa}
\usepackage{graphicx}
\usepackage{amsmath}
\usepackage{amsthm}
\usepackage{amsfonts}
\usepackage{amssymb}

\def\currentvolume{?}
\def\currentissue{?}
\def\currentyear{200?}
\def\month{?}
\def\ppages{1--25}

\newtheorem{theorem}{Theorem}[section]
\newtheorem{proposition}[theorem]{Proposition}
\newtheorem{corollary}[theorem]{Corollary}
\newtheorem{lemma}[theorem]{Lemma}
\newtheorem{definition}[theorem]{Definition}
\renewcommand{\theequation}{\thesection.\arabic{equation}}

\DeclareMathOperator{\rge}{rge}

\title[Global bifurcation of homoclinic solutions of Hamiltonian systems]
{Global bifurcation of homoclinic solutions of Hamiltonian systems}
\author[S.~Secchi \and C.~A.~Stuart]{}

\thanks{The first author is supported by MURST, national project
\emph{Variational methods and nonlinear differential equations} and Partially supported by 
EPFL}
\email{secchi@sissa.it}
\email{charles.stuart@epfl.ch}

\subjclass{34C37, 34C23}
\keywords{Homoclinic orbits, bifurcation theory, Fredholm operators.}

\begin{document}
\maketitle

\setcounter{page}{1}

\centerline{\scshape S.~Secchi${}^1$ \and C.~A.~Stuart${}^2$}
\medskip

{\footnotesize
\centerline{ ${}^1$SISSA }
\centerline{ via Beirut 2/4 }
\centerline{ I-34014 Trieste, Italy }
}
\medskip

%\centerline{\scshape C.~A.~Stuart }
%\medskip

{\footnotesize
\centerline{ ${}^2$D\'epartement de Math\'ematiques }
\centerline{ \'Ecole Polytechnique F\'ed\'erale }
\centerline{ 1015 Lausanne, Switzerland }
}
\medskip

%\centerline{(Communicated by  Glenn Webb)}

\bigskip
\begin{quote}{\normalfont\fontsize{8}{10}\selectfont
{\bfseries Abstract.}
We provide global bifurcation results for a class of nonlinear Hamiltonian systems.
\par}
\end{quote}

%\newtheorem{theorem}{Theorem}[section]
%\newtheorem{acknowledgement}[theorem]{Acknowledgement}
%\newtheorem{algorithm}[theorem]{Algorithm}
%\newtheorem{axiom}[theorem]{Axiom}
%\newtheorem{case}[theorem]{Case}
%\newtheorem{claim}[theorem]{Claim}
%\newtheorem{conclusion}[theorem]{Conclusion}
%\newtheorem{condition}[theorem]{Condition}
%\newtheorem{conjecture}[theorem]{Conjecture}
%\newtheorem{corollary}[theorem]{Corollary}
%\newtheorem{criterion}[theorem]{Criterion}
%\newtheorem{definition}[theorem]{Definition}
%\newtheorem{example}[theorem]{Example}
%\newtheorem{exercise}[theorem]{Exercise}
%\newtheorem{lemma}[theorem]{Lemma}
%\newtheorem{notation}[theorem]{Notation}
%\newtheorem{problem}[theorem]{Problem}
%newtheorem{proposition}[theorem]{Proposition}%
%newtheorem{remark}[theorem]{Remark}
%newtheorem{solution}[theorem]{Solution}
%newtheorem{summary}[theorem]{Summary}
%newcommand{\e}{\varepsilon}
%numberwithin{equation}{section}
%renewcommand{\H}{H^1}
%renewcommand{\L}{L^2}
%newcommand{\R}{\mathbb{R}}
%newcommand{\Rn}{\mathbb{R}^{2N}}

%begin{document}

%\title{Global bifurcation of homoclinic solutions of Hamiltonian systems}
%\author{S.~Secchi\thanks{Supported by MURST, national project
%\emph{Variational methods and nonlinear differential equations}. Partially supported by EPFL.} 
%\\ SISSA, via Beirut 2/4\\34014 Trieste (Italy)\and C.~A.~Stuart\\D\'epartement de 
%Math\'ematiques\\\'Ecole Polytechnique F\'ed\'erale\\Lausanne (CH)}
%\date{}
%\maketitle

\section{Introduction}

In this paper we present some results about the bifurcation of global branches
of homoclinic solutions for the following class of Hamiltonian systems:
\begin{equation}
Jx^{\prime}(t)=\nabla H(t,x(t),\lambda),\label{1.1}%
\end{equation}
where $x\in H^{1}(\mathbb{R},\mathbb{R}^{2N})$, $J$ is a real $2N\times2N-$
matrix such that $J^{T}=J^{-1}=-J$ and the Hamiltonian $H\colon\mathbb{R}%
\times\mathbb{R}^{2N}\times\mathbb{R}\longrightarrow\mathbb{R}$ is
sufficiently smooth. Moreover $\lambda$ is the bifurcation parameter and
$\nabla H(t,\xi,\lambda)=D_{\xi}H(t,\xi,\lambda)$ for $t\in\mathbb{R}$,
$\xi\in\mathbb{R}^{2N}$ and $\lambda\in\mathbb{R}$. We suppose that $x\equiv0$
satisfies (\ref{1.1}) for all values of the real parameter $\lambda$ and we
study the existence of solutions which are homoclinic to this trivial solution
in the sense that
\begin{equation}
\lim_{t\to -\infty}x(t)=\lim_{t\to +\infty}x(t)=0.\label{1.2}%
\end{equation}

Our approach is based on the topological degree for proper Fredholm operators
of index zero, as developed by Fitzpatrick, Pejsachowicz and Rabier in
\cite{fp,fpr,pr}. This tool has been applied recently by Rabier and Stuart
(see \cite{rs1}) to get bifurcation results for some classes of quasilinear
elliptic partial differential equations on $\mathbb{R}^{N}$ with possibly a
non-variational structure.

A first step is to express the problem (\ref{1.1})(\ref{1.2}) as the set of
zeros of some suitable function $F\in C^{1}(\mathbb{R}\times X,Y)$ where $X$
and $Y$ are real Banach spaces. Then we have to find conditions on $H$ under
which this $F$ is a proper Fredholm operator of index zero. The Fredholm
property holds provided that the linearisation of (\ref{1.1}) at the $x=0$
tends to periodic linear systems
\[
Ju^{\prime}(t)=A_{\lambda}^{\pm}(t)u(t)\text{ as }t\rightarrow\pm\infty
\]
which have no characteristic multipliers on the unit circle. This is proved in
Theorem \ref{th5.1}. A criterion for properness is obtained provided that the
nonlinear system (\ref{1.1}) tends to periodic (possibly autonomous)
Hamiltonian systems as $t\rightarrow\pm\infty$ which have no solutions
homoclinic to zero. The precise statement of this result is given as Theorem
\ref{th4.9} and conditions which can be used to check for the absence of
homoclinics are established in Theorem \ref{th5.5}. \ A general result
concerning the global bifurcation of solutions of the system (\ref{1.1}%
)(\ref{1.2}) is then formulated as Theorem \ref{th5.5} and we give one example
illustrating how its hypotheses can be checked.

\bigskip

\noindent\textbf{Notation}

\begin{itemize}
\item $B(X,Y)$ is the space of bounded linear operators from $X$ into $Y$.

\item $GL(X,Y)$ is the space of continuous isomorphisms from $X$ into $Y.$

\item $\deg$ is the ordinary Leray--Schauder degree.

\item  The kernel of a linear operator $L$ is denoted by $\ker L$, and its
range by $\rge L$.

\item  An operator $L\in B(X,Y)$ is said to be \emph{Fredholm of index zero}
if $\rge L$ is closed in $Y$, $\ker L$ is finite-dimensional and $\dim\ker
L=\operatorname{codim}$ $\rge L$. We set $\Phi_{0}(X,Y)=\{L\in B(X,Y):L$ is a
Fredholm operator of index zero$\}.$

\item $L^{2}=L^{2}(\mathbb{R},\mathbb{R}^{2N})$ with $\left\|  x\right\|
_{2}=(\int_{\mathbb{R}}\left\|  x(t)\right\|  ^{2}dt )^{1/2}$ for $x\in
L^{2}$ where $\left\|  \cdot\right\|  $ denotes the Euclidean norm on
$\mathbb{R}^{2N}.$ The scalar product on $\mathbb{R}^{2N}$ will be denoted by
$\left\langle \cdot,\cdot\right\rangle $ and that on $L^{2}$ by $\left\langle
\cdot,\cdot\right\rangle _{2}.$ Thus
\[
\left\langle x,y\right\rangle _{2}=\int_{-\infty}^{\infty}\left\langle
x(t),y(t)\right\rangle dt\text{ for }x,y\in L^{2}.
\]

\item $H^{1}=H^{1}(\mathbb{R},\mathbb{R}^{2N})$ with $\left\|  x\right\|
=\{\left\|  x\right\|  _{2}^{2}+\left\|  x^{\prime}\right\|  _{2}^{2}\}^{1/2}$
for $x\in H^{1}.$ Recall that, for all $x\in H^{1},$ $x$ is continuous (after
modification on a set of measure zero) and $\lim_{\left|  t\right|
\rightarrow\infty}x(t)=0.$

\item $C_{d}=\{x\in C(\mathbb{R},\mathbb{R}^{2N}):\lim_{\left|  t\right|
\rightarrow\infty}x(t)=0\}$ is Banach space with the norm $\left\|  x\right\|
_{\infty}=\sup_{t\in\mathbb{R}}\left\|  x(t)\right\|  .$ $H^{1}$ is
continuously embedded in $C_{d}.$

\item $\left\|  M\right\|  $ will also be used to denote the Euclidean norm of
a matrix $M.$
\end{itemize}

\section{A review of the topological degree for Fredholm maps}

Consider two real Banach spaces $X$ and $Y$. The notion of topological degree
for $C^{1}-$Fredholm operator of index zero from $X$ to $Y$ has been
introduced in \cite{fp,fpr,pr} in several steps. First of all, one defines the
\textbf{parity} of a continuous path $\lambda\in\lbrack a,b]\mapsto
A(\lambda)$ of bounded linear Fredholm operators with index zero from $X$ into
$Y$. It is always possible to find a \textbf{parametrix} for this path, namely
a continuous function $B\colon\lbrack a,b]\rightarrow GL(Y,X)$ such that the
composition $B(\lambda)A(\lambda)\colon X\rightarrow X$ is a compact
perturbation of the identity for every $\lambda\in\lbrack a,b]$. If $A(a)$ and
$A(b)$ belong to $GL(X,Y)$, then the parity of the path $A$ on $[a,b]$ is by
definition
\[
\pi(A(\lambda)\mid\lambda\in\lbrack a,b])=\deg(B(a)A(a))\deg(B(b)A(b)).
\]
This is a good definition in the sense that it is independent of the
parametrix $B$. The following criterion can be useful for evaluating the
parity of an admissible path.

\begin{proposition}
\label{prop2.1} Let $A\colon\lbrack a,b]\rightarrow B(X,Y)$ be a continuous
path of bounded linear operators having the following properties.

\begin{enumerate}
\item [(i)]$A\in C^{1}([a,b],B(X,Y))$.

\item[(ii)] $A(\lambda)\colon X\rightarrow Y$ is a Fredholm operator of index
zero for each $\lambda\in\lbrack a,b]$.

\item[(iii)] There exists $\lambda_{0}\in(a,b)$ such that
\begin{equation}
A^{\prime}(\lambda_{0})[\ker A(\lambda_{0})]\oplus\rge A(\lambda
_{0})=Y\label{2.1}%
\end{equation}
in the sense of a topological direct sum.
\end{enumerate}

Then there exists $\varepsilon>0$ such that $[\lambda_{0}-\varepsilon
,\lambda_{0}+\varepsilon]\subset\lbrack a,b]$,
\begin{equation}
A(\lambda)\in GL(X,Y)\text{ for }\lambda\in\lbrack\lambda_{0}-\varepsilon
,\lambda_{0})\cup(\lambda_{0},\lambda_{0}+\varepsilon]\label{2.2}%
\end{equation}
and
\begin{equation}
\pi(A(\lambda)\mid\lambda\in\lbrack\lambda_{0}-\varepsilon,\lambda
_{0}+\varepsilon])=(-1)^{k}\label{2.3}%
\end{equation}
where $k=\dim\ker A(\lambda_{0})$.
\end{proposition}

\noindent The proof of this proposition is essentially contained in
\cite{fitz,fp}.

We remark that given a continuous path $A\colon\lbrack a,b]\rightarrow\Phi
_{0}(X,Y)$ and any $\lambda_{0}\in\lbrack a,b]$ such that $A(\lambda)\in
GL(X,Y)$ for all $\lambda\neq\lambda_{0}$, the parity $\pi(A(\lambda
)\mid\lambda\in\lbrack\lambda_{0}-\varepsilon,\lambda_{0}+\varepsilon])$ is
the same for all $\varepsilon>0$ sufficiently small. This number is then
called the parity of $A$ across $\lambda_{0}$.

As in the case of the Leray--Schauder degree, the parity plays a role in
bifurcation theory.

\begin{definition}
Let $X$ and $Y$ be real Banach spaces and consider a function $F\in
C^{1}(\Lambda\times X,Y)$ where $\Lambda$ is an open interval. Let
$P(\lambda,x)=\lambda$ be the projection of $\mathbb{R}\times X$ onto
$\mathbb{R}$. We say that $\Lambda$ is an admissible interval for $F$ provided that

\begin{enumerate}
\item [(i)]for all $(\lambda,x)\in\Lambda\times X$, the bounded linear
operator $D_{x}F(\lambda,x)\colon X\rightarrow Y$ is a Fredholm operator of
index zero;

\item[(ii)] for any compact subset $K\subset Y$ and any closed bounded subset
$W$ of $\mathbb{R}\times X$ such that
\[
\inf\Lambda<\inf PW\leq\sup PW<\sup\Lambda,
\]
$F^{-1}(K)\cap W$ is a compact subset of $\mathbb{R}\times X$.
\end{enumerate}
\end{definition}

\begin{theorem}
\label{th2.3} Let $X$ and $Y$ be real Banach spaces and consider a function
$F\in C^{1}(\Lambda\times X,Y)$ where $\Lambda$ is an admissible open interval
for $F$. Suppose that $\lambda_{0}\in\Lambda$ and that there exists
$\varepsilon>0$ such that $[\lambda_{0}-\varepsilon,\lambda_{0}+\varepsilon
]\subset\Lambda$,
\[
D_{x}F(\lambda,0)\in GL(X,Y)\text{ for }\lambda\in\lbrack\lambda
_{0}-\varepsilon,\lambda_{0}+\varepsilon]\setminus\{\lambda_{0}\}
\]
and
\[
\pi(D_{x}F(\lambda,0)\mid\lbrack\lambda_{0}-\varepsilon,\lambda_{0}%
+\varepsilon])=-1.
\]
Let $Z=\{(\lambda,u)\in\Lambda\times X\mid u\neq0$ and $F(\lambda,u)=0\}$ and
let $C$ denote the connected component of $Z\cup\{(\lambda_{0},0)\}$
containing $(\lambda_{0},0)$.

Then $C$ has at least one of the following properties:

\begin{enumerate}
\item [(1)]$C$ is unbounded.

\item[(2)] The closure of $C$ contains a point $(\lambda_{1},0)$ where
$\lambda_{1}\in\Lambda\setminus\lbrack\lambda_{0}-\varepsilon,\lambda
_{0}+\varepsilon]$ and $D_{x}F(\lambda_{1},0)\notin GL(X,Y)$.

\item[(3)] The closure of $PC$ intersects the boundary of $\Lambda$.
\end{enumerate}
\end{theorem}

\begin{proof}
See \cite{rs1}.
\end{proof}

In the rest of the present paper, we want to present some explicit conditions
under which Theorem \ref{th2.3} can be applied to our problem (\ref{1.1}%
)(\ref{1.2}).

\section{The functional setting}

In our analysis of the problem, the following terminology will help us to
formulate conditions on the Hamiltonian $H$ which ensure that the system
(\ref{1.1})(\ref{1.2}) is equivalent to an equation of the form $F(\lambda,x)=0$
where $F\in C^{1}(\mathbb{R}\times H^{1},L^{2}).$

Consider a function $f\colon\mathbb{R}\times\mathbb{R}^{M}\rightarrow
\mathbb{R}.$ This $f$ can be identified with the application
\[
(t,\xi)\in\mathbb{R}\times\mathbb{R}^{M}\mapsto(t,f(t,\xi))\in\mathbb{R}%
\times\mathbb{R}.
\]

\begin{definition}
We say that $f$ is an equicontinuous $C_{\xi}^{0}$--bundle map if $f$ is
continuous on $\mathbb{R}\times\mathbb{R}^{M}$ and the collection
$\{f(t,\cdot)\}_{t\in\mathbb{R}}$ is equicontinuous at every point $\xi$ of
$\mathbb{R}^{M}$. For $k\in\mathbb{N}$, we say that $f$ is an equicontinuous
$C_{\xi}^{k}$--bundle map if all the partial derivatives $\partial
^{\kappa}f/ \partial\xi^{\kappa}$ exist for all $|\kappa|\leq k$ and are
equicontinuous $C_{\xi}^{0}$--bundle maps.
\end{definition}

We shall discuss the system (\ref{1.1})(\ref{1.2}) under the following
hypotheses on the Hamiltonian $H(t,\xi,\lambda)$ where $t,\lambda\in
\mathbb{R}$ and $\xi\in\mathbb{R}^{2N}.$

\begin{itemize}
\item [(H1)]$H\in C(\mathbb{R\times R}^{2N}\times\mathbb{R})$ with
$H(t,\cdot,\lambda)\in C^{2}(\mathbb{R}^{2N})$ and $D_{\xi}H(t,0,\lambda)=0$
for all $t,\lambda\in\mathbb{R}$

\item[(H2)] The partial derivatives $D_{\xi}H,D_{\xi}^{2}H,D_{\lambda}D_{\xi
}H,D_{\lambda}D_{\xi}^{2}H$ and $D_{\xi}D_{\lambda}D_{\xi}H$ exist and are
continuous on $\mathbb{R\times R}^{2N}\times\mathbb{R}$.

\item[(H3)] For each $\lambda\in\mathbb{R},$ $D_{\xi}H(\cdot,\cdot
,\lambda):\mathbb{R}\times\mathbb{R}^{2N}\rightarrow\mathbb{R}^{2N}$ is a
$C_{\xi}^{1}-$bundle map, and $D_{\lambda}D_{\xi}^{2}H:\mathbb{R\times
}(\mathbb{R}^{2N}\times\mathbb{R})\mathbb{\rightarrow R}$ is a $C_{(\xi
,\lambda)}^{0}-$bundle map.

\item[(H4)] $D_{\xi}^{2}H(\cdot,0,0)$ and $D_{\lambda}D_{\xi}^{2}%
H(\cdot,0,0)\in L^{\infty}(\mathbb{R}).$
\end{itemize}

\bigskip

Under these hypotheses the system (\ref{1.1}) can be expressed as
$F(\lambda,x)=0$ where
\[
F(\lambda,x)=Jx^{\prime}-h(\lambda,x)
\]
and
\[
h(\lambda,x)(t)=D_{\xi}H(t,x(t),\lambda)\text{ for }t\in\mathbb{R}%
\]
is the Nemytskii operator generated by the function $D_{\xi}H.$ To proceed we
must establish some basic properties of this Nemytskii operator. First we
observe that
\[
D_{\xi}H(t,\xi,\lambda)=\int_{0}^{1}\frac{d}{ds}D_{\xi}H(t,s\xi,\lambda
)ds=\int_{0}^{1}D_{\xi}^{2}H(t,s\xi,\lambda)\xi ds,
\]
and so
\begin{equation}
\left\|  D_{\xi}H(t,\xi,\lambda)\right\|  \leq\left\|  \xi\right\|  \int
_{0}^{1}\left\|  D_{\xi}^{2}H(t,s\xi,\lambda)\right\|  ds.\label{3.1}%
\end{equation}

\begin{lemma}
\label{lem3.2} (i) For every $K>0,$ there is a constant $C(K)>0$ such that
\[
\left\|  D_{\lambda}D_{\xi}^{2}H(t,\xi,\lambda)\right\|  \leq C(K)\text{ for
all }t\in\mathbb{R}\text{ and }(\xi,\lambda)\in\mathbb{R}^{2N+1}\text{ with
}\left\|  (\xi,\lambda)\right\|  \leq K.
\]
(ii) For each $\lambda\in\mathbb{R}$ and $K>0,$ there exists a constant
$C(\lambda,K)$ such that
\[
\left\|  D_{\xi}^{2}H(t,\xi,\lambda)\right\|  \leq C(\lambda,K)\text{ for all
}t\in\mathbb{R}\text{ and }\xi\in\mathbb{R}^{2N}\text{ with }\left\|
\xi\right\|  \leq K.
\]
\end{lemma}

\begin{proof}(i) The hypothesis (H4) means that$\ \sup
_{t\in\mathbb{R}} \|  D_{\lambda}D_{\xi}^{2}H(t,0,0) \|  <\infty,$
and (H3) implies that, for all $(\xi,\lambda)\in\mathbb{R}^{2N+1},$ there
exists $\delta(\xi,\lambda)>0$ such that
\[
\left\|  D_{\lambda}D_{\xi}^{2}H(t,\xi,\lambda)-D_{\lambda}D_{\xi}^{2}%
H(t,\eta,\mu)\right\|  <1\text{ for all }t\in\mathbb{R}%
\]
provided that $\left\|  (\xi,\lambda)-(\eta,\mu)\right\|  <\delta(\xi
,\lambda).$ A straight forward compactness argument now leads to the first assertion.

(ii) First we note that
\[
D_{\xi}^{2}H(t,0,\lambda)-D_{\xi}^{2}H(t,0,0)=\int_{0}^{1}\frac{d}{ds}D_{\xi
}^{2}H(t,0,s\lambda)ds=\lambda\int_{0}^{1}D_{\lambda}D_{\xi}^{2}%
H(t,0,s\lambda)ds.
\]
Hence, from (H4) and part (i), we see that
\[
\sup_{t\in\mathbb{R}}\left\|  D_{\xi}^{2}H(t,0,\lambda)\right\|  \leq
\sup_{t\in\mathbb{R}}\left\|  D_{\xi}^{2}H(t,0,0)\right\|  +\left|
\lambda\right|  C(\left|  \lambda\right|  ).
\]
Furthermore, by (H3), for any $\xi\in\mathbb{R}^{2N},$ there exists
$\delta(\xi,\lambda)>0$ such that
\[
\left\|  D_{\xi}^{2}H(t,\xi,\lambda)-D_{\xi}^{2}H(t,\eta,\lambda)\right\|
<1\text{ for all }t\in\mathbb{R}%
\]
provided that $\left\|  \xi-\eta\right\|  <\delta(\xi,\lambda).$ A compactness
argument now yields the conclusion (ii).

This lemma shows that the Nemytskii operator $h(\lambda,\cdot)$ maps $H^{1} $
into $L^{2}.$ Indeed, for any $x\in H^{1},$ we have that $\left\|  x\right\|
_{\infty}<\infty,$ and so by the lemma, there exists a constant $C(\lambda
,\left\|  x\right\|  _{\infty})$ such that
\[
\left\|  D_{\xi}^{2}H(t,\xi,\lambda)\right\|  \leq C(\lambda,\left\|
x\right\|  _{\infty})\text{ for all }t\in\mathbb{R}\text{ and }\left\|
\xi\right\|  \leq\left\|  x\right\|  _{\infty}%
\]
Hence, by (\ref{3.1}),
\[
\left\|  D_{\xi}H(t,x(t),\lambda)\right\|  \leq\left\|  x(t)\right\|  \int
_{0}^{1}\left\|  D_{\xi}^{2}H(t,sx(t),\lambda)\right\|  ds\leq\left\|
x(t)\right\|  C(\lambda,\left\|  x\right\|  _{\infty}),
\]
showing that
\begin{equation}
h(\lambda,x)\in L^{2}\text{ with }\left\|  h(\lambda,x)\right\|  _{2}\leq
C(\lambda,\left\|  x\right\|  _{\infty})\left\|  x(t)\right\|  _{2}%
.\label{3.2}%
\end{equation}

From now on we can consider $h$ as a mapping from $\mathbb{R}\times H^{1}$
into $L^{2}$ and the system (\ref{1.1})(\ref{1.2}) can be written as
\begin{align}
F(\lambda,x) &  =0\text{ where }F:\mathbb{R}\times H^{1}\rightarrow
L^{2}\text{ is defined by}\nonumber\\
F(\lambda,x) &  =Jx^{\prime}-h(\lambda,x).\label{3.3}%
\end{align}
Note that if $(\lambda,x)\in\mathbb{R}\times H^{1}$ and $F(\lambda,x)=0,$ it
follows that $x\in C^{1}(\mathbb{R})$ and $\lim_{\left|  t\right|
\rightarrow\infty}x(t)=0.$ Furthermore, it follows from the argument leading
to (\ref{3.2}) that $F(\lambda,\cdot)$ maps $H^{1}$ boundedly into $L^{2}.$
\end{proof}

We now investigate the smoothness of the function $F:\mathbb{R}\times
H^{1}\rightarrow L^{2}.$ The Hamiltonian $H$ is assumed to have the properties
(H1) to (H4) from now on.

\begin{theorem}
\label{th3.3}

\begin{itemize}
\item [(1)]$F\in C^{1}(\mathbb{R}\times H^{1},L^{2})$ with
\[
D_{x}F(\lambda,x)u=Ju^{\prime}-M(\lambda,x)u\text{ for all }\lambda
\in\mathbb{R}\text{ and }x,u\in H^{1}%
\]
where $M(\lambda,x)(t)=D_{\xi}^{2}H(t,x(t),\lambda)$ for all $t\in\mathbb{R}.$

\item[(2)] $D_{x}F(\cdot,0)\in C^{1}(\mathbb{R},B(H^{1},L^{2}))$ and
\[
D_{\lambda}D_{x}F(\lambda,0)u=-D_{\lambda}D_{\xi}^{2}H(t,0,\lambda)u\text{ for
all }\lambda\in\mathbb{R}\text{ and }u\in H^{1}.
\]

\item[(3)] Let $W$ be any bounded subset of $H^{1}.$ The family of functions
$\{F(\cdot,u):\mathbb{R}\rightarrow L^{2}\}_{u\in W}$ is equicontinuous at
$\lambda$ for every $\lambda\in\mathbb{R}.$

\item[(4)] For all $\lambda\in\mathbb{R},$ the function $F(\lambda
,\cdot):H^{1}\rightarrow L^{2}$ is weakly sequentially continuous.
\end{itemize}
\end{theorem}

\begin{proof}See the Appendix.
\end{proof}

Noting that $M(\lambda,x)(t)$ is a symmetric $2N\times2N-$matrix, we see that
the equation
\[
D_{x}F(\lambda,x)u=0
\]
is a linear Hamiltonian system. We have already established that, for all
$\lambda\in\mathbb{R}$ and $x\in H^{1},$ $D_{x}F(\lambda,x):H^{1}\rightarrow
L^{2}$ is a bounded linear operator. It is important to know when it is a
Fredholm operator of index zero. In fact, $D_{x}F(\lambda,x)$ can also be
considered as an unbounded self-adjoint operator acting in $L^{2}$ and this
means that its index must be zero whenever it is Fredholm. The next result
summarizes the situation. Later we shall give explicit conditions on $H$ which
ensure that $D_{x}F(\lambda,x)\in\Phi_{0}(H^{1},L^{2}).$

\begin{theorem}
\label{th3.4} For $\lambda\in\mathbb{R}$ and $x\in H^{1},$ set $L=D_{x}%
F(\lambda,x)$.

\begin{itemize}
\item [(1)]$L\in B(H^{1},L^{2})$

\item[(2)] $L:H^{1}\subset L^{2}\rightarrow L^{2}$ is an unbounded
self-adjoint operator in $L^{2}$, and so
\[
L^{2}=\ker L\oplus\overline{\rge L}%
\]
where $\overline{\rge L}$ is the closure of the range of $L$ in $L^{2}$
and $\oplus$ denotes an orthogonal direct sum in $L^{2}.$

\item[(3)] The following statements are equivalent:

\begin{itemize}
\item [(a)]$\rge L=\overline{\rge L}$

\item[(b)] $L\in\Phi_{0}(H^{1},L^{2}).$
\end{itemize}

\item[(4)] Let $w\in H^{1}.$ Then $L=D_{x}F(\lambda,x)\in\Phi_{0}(H^{1}%
,L^{2})\Longleftrightarrow D_{x}F(\lambda,w)\in\Phi_{0}(H^{1},L^{2}).$
\end{itemize}
\end{theorem}

\begin{proof}(1) For any $u\in H^{1},$%
\[
\left\|  Lu\right\|  _{2}=\left\|  Ju^{\prime}-M(\lambda,x)u\right\|  _{2}%
\leq\left\|  Ju^{\prime}\right\|  _{2}+\left\|  M(\lambda,x)u\right\|  _{2}%
\]
by Theorem \ref{th3.3} where $\left\|  M(\lambda,x)(t)\right\|  \leq
C(\lambda,\left\|  x\right\|  _{\infty})$ for all $t\in\mathbb{R}$ by Lemma
\ref{lem3.2} Thus
\[
\left\|  Lu\right\|  _{2}\leq\left\|  u^{\prime}\right\|  _{2}+C(\lambda
,\left\|  x\right\|  _{\infty})\left\|  u\right\|  _{2}%
\]
showing that $L\in B(H^{1},L^{2}).$

(2) This is well-known. See \cite{joost}, for example.

(3) Clearly, (b)$\Longrightarrow$(a) and so it suffices to prove that
(a)$\Longrightarrow$(b). First we observe that if $u\in\ker L,$ the $u\in
C^{1}(\mathbb{R})$ and $Ju^{\prime}=M(\lambda,x)u.$ But the set of all
solutions of this linear system is a vector space of dimension $2N$ and hence
$\dim\ker L\leq2N.$ Furthermore, by part (2), $\operatorname{codim}%
\overline{\rge L}=\dim\ker L.$ Thus (a)$\Longrightarrow$(b).

(4) For all $x,w,u\in H^{1},$%
\[
\{D_{x}F(\lambda,x)-D_{x}F(\lambda,w)\}u=\{M(\lambda,w)-M(\lambda,x)\}u
\]
where, for all $t\in\mathbb{R},$%
\begin{multline*}
M(\lambda,w)(t)-M(\lambda,x)(t) =\\
\{D_{\xi}^{2}H(t,w(t),\lambda)-D_{\xi}^{2}H(t,0,\lambda)\}-\{D_{\xi}%
^{2}H(t,x(t),\lambda)-D_{\xi}^{2}H(t,0,\lambda)\}.
\end{multline*}
Given any $\varepsilon>0,$ it follows from (H3) that there exists $\delta>0 $
such that
\[
\left\|  D_{\xi}^{2}H(t,\xi,\lambda)-D_{\xi}^{2}H(t,0,\lambda)\right\|
<\varepsilon
\]
for all $t\in\mathbb{R}$ and for all $\xi\in\mathbb{R}^{2N}$ such that
$\left|  \xi\right|  <\delta.$ But, since $x,w\in H^{1},$ there exists $R>0$
such that $\left|  x(t)\right|  <\delta$ and $\left|  w(t)\right|  <\delta$
whenever $\left|  t\right|  >R.$ Thus,
\[
\left\|  M(\lambda,w)(t)-M(\lambda,x)(t)\right\|  <2\varepsilon\text{ for all
}t\in\mathbb{R}\text{ such that }\left|  t\right|  >R,
\]
showing that $\lim_{\left|  t\right|  \rightarrow\infty}\left\|
M(\lambda,w)(t)-M(\lambda,x)(t)\right\|  =0.$ From this it follows easily that
multiplication by $M(\lambda,w)-M(\lambda,x)$ defines a compact linear
operator $K$ from $H^{1}$ into $L^{2}.$ Since $D_{x}F(\lambda,x)-D_{x}%
F(\lambda,w)=K$, this implies that $D_{x}F(\lambda,x)\in\Phi_{0}(H^{1}%
,L^{2})\Longleftrightarrow D_{x}F(\lambda,w)\in\Phi_{0}(H^{1},L^{2})$.
\end{proof}

\section{Admissible intervals}

In this section we give some useful criteria for the existence of admissible
intervals for $F.$ For this we shall assume henceforth that, in addition to
the properties (H1) to (H4), the Hamiltonian $H$ is asymptotically periodic in
the following sense.

(H$^{\infty}$) For all $\lambda\in\mathbb{R},$ there exist two $C_{\xi}^{1}%
-$bundle maps $g^{+}(\cdot,\cdot,\lambda)$ and $g^{-}(\cdot,\cdot
,\lambda):\mathbb{R}\times\mathbb{R}^{2N}\mathbb{\rightarrow R}^{2N}$ such that

\begin{itemize}
\item [(1)]$g^{+}(t,0,\lambda)=g^{-}(t,0,\lambda)=0\text{ for all }%
t,\lambda\in\mathbb{R}$

\medskip

\item[(2)] $\lim_{t\rightarrow\infty}\{D_{\xi}^{2}H(t,\xi,\lambda)-D_{\xi
}g^{+}(t,\xi,\lambda)\}=\lim_{t\rightarrow-\infty}\{D_{\xi}^{2}H(t,\xi
,\lambda)-\newline D_{\xi}g^{-}(t,\xi,\lambda)\}=0$, uniformly for $\xi$ in
bounded subsets of $\mathbb{R}^{2N}$

\medskip

\item[(3)] $g^{+}(t+T^{+},\xi,\lambda)-g^{+}(t,\xi,\lambda)=g^{-}(t+T^{-}%
,\xi,\lambda)-g^{-}(t,\xi,\lambda)=0$ for some $T^{+},T^{-}>0$ and for all
$(t,\xi)\in\mathbb{R}\times\mathbb{R}^{2N}$.
\end{itemize}

\noindent\textbf{Remarks} (1) The periods $T^{+}$ and $T^{-}$ may depend on
$\lambda.$

(2) It follows easily from this assumption that $D_{\xi}g^{+}(t,\xi,\lambda)$
and $D_{\xi}g^{-}(t,\xi,\lambda)$ are symmetric matrices for all
$(t,\xi,\lambda)\in\mathbb{R}\times\mathbb{R}^{2N}\times\mathbb{R},$ and that
\begin{equation}
\lim_{t\rightarrow\infty}\{D_{\xi}H(t,\xi,\lambda)-g^{+}(t,\xi,\lambda
)\}=\lim_{t\rightarrow-\infty}\{D_{\xi}H(t,\xi,\lambda)-g^{-}(t,\xi
,\lambda)\}=0,\label{4.1}%
\end{equation}
uniformly for $\xi$ in bounded subsets of $\mathbb{R}^{2N}$.

Furthermore, setting
\[
H^{\pm}(t,\xi,\lambda)=H(t,0,\lambda)+\int_{0}^{1}\left\langle g^{\pm}%
(t,s\xi,\lambda),\xi\right\rangle ds,
\]
we have that
\begin{equation}
D_{\xi}H^{\pm}(t,\xi,\lambda)=g^{\pm}(t,\xi,\lambda)\label{4.2}%
\end{equation}
and
\begin{equation}
\lim_{t\rightarrow\infty}\{H(t,\xi,\lambda)-H^{+}(t,\xi,\lambda)\}=\lim
_{t\rightarrow-\infty}\{H(t,\xi,\lambda)-H^{-}(t,\xi,\lambda)\}=0\label{4.3}%
\end{equation}
uniformly for $\xi$ in bounded subsets of $\mathbb{R}^{2N}$. In particular,
the differential equations
\[
Jx^{\prime}(t)-g^{+}(t,x(t),\lambda)=0\text{ and }Jx^{\prime}(t)-g^{-}%
(t,x(t),\lambda)=0
\]
are periodic Hamiltonian systems. Let $h^{\pm}$ denote the Nemytskii operators
induced by $g^{\pm},h^{\pm}(\lambda,x)(t)=g^{\pm}(t,x(t),\lambda)$ and then
define $F^{\pm}$ by
\begin{equation}
F^{\pm}(\lambda,x)=Jx^{\prime}-h^{\pm}(\lambda,x).\label{4.4}%
\end{equation}

\begin{theorem}\label{th4.1}
Under the hypotheses (H1) to (H4) and (H$^{\infty})$, for every
$\lambda\in\mathbb{R},$ $F^{\pm}(\lambda,\cdot)$ maps $H^{1}$ boundedly into
$L^{2}.$ Furthermore, $F^{\pm}(\lambda,\cdot):H^{1}\rightarrow L^{2}$ is
weakly sequentially continuous.
\end{theorem}

\begin{proof} See the Appendix.
\end{proof}

\begin{lemma}
\label{lem4.2} Let $B$ be a bounded subset of $H^{1}$ and consider $\lambda
\in\mathbb{R}$ and $\varepsilon,L>0.$ There exists $R=R(\varepsilon
,B,L,\lambda)>0$ such that
\[
\left\|  F(\lambda,x)-F^{+}(\lambda,x)\right\|  _{L^{2}(I^{+})}<\varepsilon
\text{ and }\left\|  F(\lambda,x)-F^{-}(\lambda,x)\right\|  _{L^{2}(I^{-}%
)}<\varepsilon
\]
for all $x\in B$ where $I^{\pm}$ are any intervals of length less than $L$
with $I^{+}\subset\lbrack R,\infty)$ and $I^{-}\subset(-\infty,-R].$
\end{lemma}

\begin{proof}
Since $B$ is bounded in $H^{1},$ there is a constant
$b>0$ such that
\[
\left\|  x(t)\right\|  \leq b\text{ for all }t\in\mathbb{R}\text{ and all
}x\in B.
\]
By (\ref{4.1}), there exists $R=R(b,\varepsilon,\lambda,L)>0$ such that
\[
\left\|  D_{\xi}H(t,\xi,\lambda)-g^{+}(t,\xi,\lambda)\right\|  <\sqrt
{\frac{\varepsilon}{L}}\text{ for all }t\geq R\text{ and }\left\|
\xi\right\|  \leq b,
\]
and hence,
\[
\left\|  D_{\xi}H(t,x(t),\lambda)-g^{+}(t,x(t),\lambda)\right\|  <\sqrt
{\frac{\varepsilon}{L}}\text{ for all }t\geq R\text{ and }x\in B.
\]
It follows that, for any interval $I^{+}$ of length less than $L$ with
$I^{+}\subset\lbrack R,\infty),$%
\begin{align*}
\left\|  F(\lambda,x)-F^{+}(\lambda,x)\right\|  _{L^{2}(I^{+})}^{2}  &
=\int_{I^{+}}\left\|  D_{\xi}H(t,x(t),\lambda)-g^{+}(t,x(t),\lambda)\right\|
^{2}dt\\
& \leq\varepsilon\text{ for all }x\in B.
\end{align*}
The other case is similar.
\end{proof}

We now introduce a notation for the translate of a function. Given
$h\in\mathbb{R}$ and $f:\mathbb{R}\rightarrow\mathbb{R}^{M},$ let $\tau
_{h}(f)$ be the function defined by
\[
\tau_{h}(f)(t)=f(t+h)\text{ for all }t\in\mathbb{R}.
\]
In particular, $\tau_{h}(F(\lambda,x))$ is the function
\[
\tau_{h}(F(\lambda,x))(t)=Jx^{\prime}(t+h)-D_{\xi}H(t+h,x(t+h),\lambda).
\]

\begin{lemma}
\label{lem4.3}Let $B$ be a bounded subset of $H^{1}$ and consider $\lambda
\in\mathbb{R}$ and $\varepsilon,\omega>0.$ There exists $h_{0}=h_{0}%
(\varepsilon,B,\omega,\lambda)>0$ such that
\begin{align*}
\left\|  \tau_{h}(F(\lambda,x))-\tau_{h}(F^{+}(\lambda,x))\right\|
_{L^{2}(-\omega,\omega)}<\varepsilon\text{ for all }x\in B\text{ and }h\geq
h_{0}\\
\text{and }\left\|  \tau_{h}(F(\lambda,x))-\tau_{h}(F^{-}(\lambda,x))\right\|
_{L^{2}(-\omega,\omega)}<\varepsilon\text{ for all }x\in B\text{ and }%
h\leq-h_{0}.
\end{align*}
\end{lemma}

\begin{proof}
Since 
\[
\left\|  \tau_{h}(F(\lambda,x))-\tau
_{h}(F^{\pm}(\lambda,x))\right\|  _{L^{2}(-\omega,\omega)}=\left\|
F(\lambda,x)-F^{\pm}(\lambda,x)\right\|  _{L^{2}(-\omega+h,\omega+h)},
\]
 the
result follows immediately from Lemma \ref{lem4.2}.
\end{proof}

\begin{lemma}
\label{lem4.4} Let $\{h_{n}\}\subset\mathbb{R}$ be a sequence such that
$\lim_{n\rightarrow\infty}\left|  h_{n}\right|  =\infty.$ For any $x\in
L^{2},$ let $\widetilde{x_{n}}=\tau_{h_{n}}(x).$ Then$\ \tilde{x}%
_{n}\rightharpoonup0$ weakly in $L^{2}.$
\end{lemma}

\begin{proof}
Since $\left\|  \widetilde{y_{n}}\right\|
_{2}=\left\|  y\right\|  _{2}$ for all $n,$ it is enough to show that
\[
\left\langle \widetilde{y_{n}},\varphi\right\rangle _{2}\rightarrow0\text{ as
}n\rightarrow\infty\text{ for all }\varphi\in C_{0}^{\infty}(\mathbb{R}%
,\mathbb{R}^{2N})
\]
where $\left\langle \cdot,\cdot\right\rangle _{2}$ denotes the usual scalar
product on $L^{2}.$ Suppose that $\varphi(t)=0$ for all $\left|  t\right|
>L.$ Then,
\begin{align*}
\left|  \left\langle \widetilde{y_{n}},\varphi\right\rangle _{2}\right|   &
\leq\int_{-L}^{L}\left\|  \widetilde{y_{n}}(t)\right\|  \left\|
\varphi(t)\right\|  dt\leq\sqrt{2L}\left\|  \varphi\right\|  _{\infty}%
\left( \int_{-L}^{L}\left\|  \widetilde{y_{n}}(t)\right\|  ^{2}dt\right)^{1/2}\\
& =\sqrt{2L}\left\|  \varphi\right\|  _{\infty}\left( \int_{-L+h_{n}}^{L+h_{n}%
}\left\|  y(t)\right\|  ^{2}dt\right)^{1/2}%
\end{align*}
and
\[
\int_{-L+h_{n}}^{L+h_{n}}\left\|  y(t)\right\|  ^{2}dt\rightarrow0\text{ as
}n\rightarrow\infty
\]
since $y\in L^{2}$ and $\left|  h_{n}\right|  \rightarrow\infty.$
\end{proof}

\begin{definition}
We say that a sequence $\{x_{n}\}$ in $H^{1}$ \textbf{vanishes uniformly at
infinity} if, for all $\varepsilon>0,$ there exists $R>0$ such that $\left\|
x_{n}(t)\right\|  \leq\varepsilon$ for all $|t|\geq R$ and all $n\in\mathbb{R}.$
\end{definition}

Recalling that $H^{1}$ is continuously embedded in $C_{d},$ we observe that
$\{x_{n}\}\subset H^{1}$ vanishes uniformly at infinity if, for all
$\varepsilon>0,$ there exist $R>0$ and $n_{0}\in\mathbb{N}$ such that
$\left\|  x_{n}(t)\right\|  \leq\varepsilon$ for all $|t|\geq R$ and all
$n\geq n_{0}.$

\begin{lemma}
\label{lem4.6}Let $\{x_{n}\}$ be a bounded sequence in $H^{1}$ and let $x\in
H^{1}.$ The following statements are equivalent.

(1) $\left\|  x_{n}-x\right\|  _{\infty}\rightarrow0$,

(2) $x_{n}\rightharpoonup x$ weakly in $H^{1}$ and $\{x_{n}\}$ vanishes
uniformly at infinity.
\end{lemma}

\begin{proof}
We begin by showing that (1) implies (2). If
$\{x_{n}\}$ does not converge weakly to $x,$ there are a number $\delta>0$, an
element $\varphi\in H^{1}$ and a subsequence $\{x_{n_{k}}\}$ such that
\[
\left|  \left\langle x_{n_{k}}-x,\varphi\right\rangle _{2}+\left\langle
x_{n_{k}}^{\prime}-x^{\prime},\varphi^{\prime}\right\rangle _{2}\right|
\geq\delta\text{ for all }k.
\]
Then, passing to a further subsequence, we can suppose that $\{x_{n_{k}}\}$
converges weakly in $H^{1}$ to some element $y.$ This implies that
$\{x_{n_{k}}\}$ converges uniformly to $y$ on any compact interval and so, by
(1), $y=x.$ This contradicts the choice of $\{x_{n_{k}}\},$ proving that
$\{x_{n}\}$ must converge weakly to $x$ in $H^{1}.$\newline Now fix
$\varepsilon>0.$ There exists $R>0$ such that $\left\|  x(t)\right\|
<\varepsilon$ for all $\left|  t\right|  \geq R,$ and there exists $n_{0}%
\in\mathbb{N}$ such that $\left\|  x_{n}-x\right\|  _{\infty}<\varepsilon$ for
all $n\geq n_{0}.$ Hence, for all $\left|  t\right|  \geq R$ and all $n\geq
n_{0},$%
\[
\left\|  x_{n}(t)\right\|  \leq\left\|  x_{n}-x\right\|  _{\infty}+\left\|
x(t)\right\|  <2\varepsilon,
\]
showing that $\{x_{n}\}$ vanishes uniformly at infinity.

Now we show that (2) implies (1). For any $\varepsilon>0,$ there exists $R>0 $
such that
\[
\left\|  x_{n}(t)\right\|  <\varepsilon\text{ for all }\left|  t\right|  \geq
R\text{ and all }n\in\mathbb{N},
\]
since $\{x_{n}\}$ vanishes uniformly at infinity. Increasing $R$ if necessary,
we also have that
\[
\left\|  x(t)\right\|  <\varepsilon\text{ for all }\left|  t\right|  \geq R
\]
since $x\in H^{1}\subset C_{d}.$ Thus,
\[
\left\|  x_{n}(t)-x(t)\right\|  <2\varepsilon\text{ for all }\left|  t\right|
\geq R\text{ and all }n.
\]
But the weak convergence of $\{x_{n}\}$ in $H^{1}$ implies that $x_{n}$
converges uniformly to $x$ on $[-R,R],$ so there exists $n_{1}\in\mathbb{N}$
such that
\[
\left\|  x_{n}(t)-x(t)\right\|  <\varepsilon\text{ for all }\left|  t\right|
\leq R\text{ and all }n\geq n_{1}.
\]
Thus we have that
\[
\left\|  x_{n}-x\right\|  _{\infty}\leq2\varepsilon\text{ for all }n\geq
n_{1},
\]
showing that $\left\|  x_{n}-x\right\|  _{\infty}\rightarrow0$.
\end{proof}

\begin{theorem}
\label{th4.7} Recalling the hypotheses (H1) to (H4), suppose that there is an
element $(\lambda,x^{0})\in\mathbb{R}\times H^{1}$ such that $D_{x}%
(\lambda,x^{0})\in\Phi_{0}(H^{1},L^{2}).$ Then the following statements are
equivalent.\newline (1) The restriction of $F(\lambda,\cdot):H^{1}\rightarrow
L^{2}$ to the closed bounded subsets of $H^{1}$ is proper.\newline (2) Every
bounded sequence $\{x_{n}\}$ in $H^{1}$ such that $\{F(\lambda,x_{n})\}$ is
convergent in $L^{2}$ contains a subsequence converging in $C_{d}.$
\end{theorem}

\begin{proof}
We show first that (1) implies (2). Indeed, let
$\{x_{n}\}$ be a bounded sequence in $H^{1}$ such that $\{F(\lambda,x_{n})\}$
converges in $L^{2}$. We want to prove that some subsequence of $\{x_{n}\}$
converges in $C_{d}.$ We know that $\left\|  x_{n}\right\|  \leq M$ for every
$n\in\mathbb{N}$ and some constant $M>0$; moreover, there is an element $y\in
L^{2}$ such that $\left\|  F(\lambda,x_{n})-y\right\|  _{2}\rightarrow0$. Let
$K=\{F(\lambda,x_{n})\}\cup\{y\}$ and $W=\overline{B(0,2M)}\subset H^{1}$.
Then $K$ is compact in $L^{2}$ and $W$ is closed and bounded in $H^{1}$. From
assumption (1) we know that
\[
F(\lambda,\cdot)^{-1}(K)\cap W\text{ is compact.}%
\]
We conclude that $\{x_{n}\}$ has a strongly convergent subsequence in $H^{1}$,
and \textit{a fortiori} in $C_{d}.$

To show that (2) implies (1), we proceed as follows. Let $W$ be a closed and
bounded subset$\ $of $H^{1}$, and let $K$ be a compact subset of $L^{2}$. We
wish to prove that
\[
F(\lambda,\cdot)^{-1}(K)\cap W\text{ is compact.}%
\]
Let $\{x_{n}\}$ be a sequence from $F(\lambda,\cdot)^{-1}(K)\cap W$. In
particular, there exists a constant $M>0$ such that $\Vert x_{n}\Vert\leq M$
for all $n$. Moreover, passing to a subsequence, we can suppose that
\[
F(\lambda,x_{n})\rightarrow y\in K\text{ strongly in }L^{2}.
\]
We can also assume that $x_{n}\rightharpoonup x$ in $H^{1},$ and by (2) that
$\left\|  x_{n}-x\right\|  _{\infty}\rightarrow0.$ By the weak sequential
continuity of $F(\lambda,\cdot)$, we get $F(\lambda,x_{n})\rightharpoonup
F(\lambda,x)$ in $L^{2}$, so that $y=F(\lambda,x)$.

We claim that
\begin{equation}
\lim_{n\rightarrow+\infty}\Vert F(\lambda,x_{n})-F(\lambda,x)-D_{x}%
F(\lambda,x)(x_{n}-x)\Vert_{L^{2}}=0.\label{4.5}%
\end{equation}
Indeed,
\begin{multline*}
F(\lambda,x_{n})-F(\lambda,x)-D_{x}F(\lambda,x)(x_{n}-x)\\
=-D_{\xi}H(\cdot,x_{n},\lambda)+D_{\xi}H(\cdot,x,\lambda)+D_{\xi}^{2}%
H(\cdot,x,\lambda)(x_{n}-x)\\
=D_{\xi}^{2}H(\cdot,x,\lambda)(x_{n}-x)-\int_{0}^{1}\frac{d}{d\tau}D_{\xi
}H(\cdot,\tau x_{n}+(1-\tau)x,\lambda)\,d\tau\\
=D_{\xi}^{2}H(\cdot,x,\lambda)(x_{n}-x)-\int_{0}^{1}D_{\xi}^{2}H(\cdot,\tau
x_{n}+(1-\tau)x,\lambda)(x_{n}-x)\,d\tau\\
=\int_{0}^{1}\{D_{\xi}^{2}H(\cdot,x,\lambda)-D_{\xi}^{2}H(\cdot,\tau
x_{n}+(1-\tau)x,\lambda)\}(x_{n}-x)\,d\tau.
\end{multline*}
Since $D_{\xi}H(\cdot,\cdot,\lambda)$ is a $C_{\xi}^{1}$--bundle map, a
standard compactness argument shows that, for any $\varepsilon>0,$ there is a
$\delta>0$ such that
\[
\sup_{t\in\mathbb{R}}\Vert D_{\xi}^{2}H(t,\xi,\lambda)-D_{\xi}^{2}%
H(t,\eta,\lambda)\Vert<\varepsilon
\]
for all $\xi,\eta\in\mathbb{R}^{2N}$ such that $\left\|  \xi\right\|
,\left\|  \eta\right\|  \leq\left\|  x\right\|  _{\infty}+1$ and $\left\|
\xi-\eta\right\|  <\delta.$ But 
\[
\left\|  [\tau x_{n}+(1-\tau)x]-x\right\|
_{\infty}=\tau\left\|  x_{n}-x\right\|  _{\infty}\leq\left\|  x_{n}-x\right\|
_{\infty}
\]
for all $\tau\in\lbrack0,1],$ and so there exists $n_{0}%
\in\mathbb{N}$ such that $\left\|  [\tau x_{n}+(1-\tau)x]-x\right\|  _{\infty
}<\delta$ for all $n\geq n_{0}$ and all $\tau\in\lbrack0,1]$. Thus,
\[
\sup_{t\in\mathbb{R}}\Vert D_{\xi}^{2}H(t,x(t),\lambda)-D_{\xi}^{2}H(t,\tau
x_{n}(t)+(1-\tau)x(t),\lambda)\Vert<\varepsilon
\]
for all $n\geq n_{0}$ and all $\tau\in\lbrack0,1]$. Hence,
\[
\left\|  F(\lambda,x_{n})-F(\lambda,x)-D_{x}F(\lambda,x)(x_{n}-x)\right\|
_{2}\leq\varepsilon\left\|  x_{n}-x\right\|  _{2}\leq\varepsilon2M
\]
for all $n\geq n_{0}$, since $\left\|  \cdot\right\|  _{2}\leq\left\|
\cdot\right\|  .$ This proves that
\[
\left\|  F(\lambda,x_{n})-F(\lambda,x)-D_{x}F(\lambda,x)(x_{n}-x)\right\|
_{2}\rightarrow0\text{ as }n\rightarrow\infty.
\]
But $\left\|  F(\lambda,x_{n})-F(\lambda,x)\right\|  _{2}=\left\|
F(\lambda,x_{n})-y\right\|  _{2}\rightarrow0$ and so we have that
\[
\left\|  L(x_{n}-x)\right\|  _{2}\rightarrow0\text{ where }L=D_{x}%
F(\lambda,x)\in\Phi_{0}(H^{1},L^{2}).
\]
By Theorem 3.15 of \cite{ee}, there exist $S\in B(L^{2},H^{1})$ and a compact
linear operator $C:H^{1}\rightarrow H^{1}$ such that $SL=I+C.$ But then,
\begin{align*}
\left\|  x_{n}-x\right\|   & =\left\|  (SL-C)(x_{n}-x)\right\|  \leq\left\|
SL(x_{n}-x)\right\|  +\left\|  C(x_{n}-x)\right\| \\
& \leq\left\|  S\right\|  \left\|  L(x_{n}-x)\right\|  _{2}+\left\|
C(x_{n}-x)\right\|
\end{align*}
where $\left\|  C(x_{n}-x)\right\|  \rightarrow0$ by the compactness of $C$
and the weak convergence of $\{x_{n}\}$ in $H^{1}.$ Thus $\left\|
x_{n}-x\right\|  \rightarrow0$ and the compactness of $F(\lambda,\cdot
)^{-1}(K)\cap W$ is established.
\end{proof}

\begin{lemma}
\label{lem4.8} Let $\{x_{n}\}$ be a bounded sequence in $H^{1}$ and consider
any numbers $T^{+},T^{-}\in(0,\infty).$ At least one of the following
properties must hold.

(1) $\{x_{n}\}$ vanishes uniformly at infinity

(2) There is a sequence $\{l_{k}\}\subset\mathbb{Z}$ with $\lim_{k\rightarrow
\infty}l_{k}=\infty$ and a subsequence $\{x_{n_{k}}\}$ of $\{x_{n}\}$ such
that $\widetilde{x_{k}}=\tau_{l_{k}T^{+}}(x_{n_{k}})=x_{n_{k}}(\cdot
+l_{k}T^{+})$ converges weakly in $H^{1}$ to an element $\widetilde{x}\neq0$

(3) There is a sequence $\{l_{k}\}\subset\mathbb{Z}$ with $\lim_{k\rightarrow
\infty}l_{k}=-\infty$ and a subsequence $\{x_{n_{k}}\}$ of $\{x_{n}\}$ such
that $\widetilde{x_{k}}=\tau_{l_{k}T^{-}}(x_{n_{k}})=x_{n_{k}}(\cdot
+l_{k}T^{-})$ converges weakly in $H^{1}$ to an element $\widetilde{x}\neq0$
\end{lemma}

\begin{proof}
Assume that $\{x_{n}\}$ does not satisfy (1). Then
there exists $\varepsilon>0$ such that, for all $k\in\mathbb{N},$ there exists
$t_{k}\in\mathbb{R}$ with $|t_{k}|\geq k$, and there exists $n_{k}%
\in\mathbb{N}$ with $n_{k}>k$ such that $\left\|  x_{n_{k}}(t_{k})\right\|
\geq\varepsilon$. By passing to a subsequence we may suppose that $\{t_{k}\}$
diverges either to $+\infty$ or to $-\infty$.\newline Suppose that
$t_{k}\rightarrow\infty$ as $k\rightarrow\infty.$ Then there exists $l_{k}%
\in\mathbb{Z}$ such that $t_{k}-l_{k}T^{+}\in\lbrack0,T^{+}]$ and
$l_{k}\rightarrow\infty$ as $k\rightarrow\infty.$ Clearly, $\Vert\tilde
{x}_{_{k}}\Vert=\Vert x_{n_{k}}\Vert$ and hence also $\{\tilde{x}_{k}\}$ is
bounded in $H^{1}.$ Passing to a further subsequence which we still denote by
$\{\tilde{x}_{k}\},$ we can suppose that $\tilde{x}_{k}$ converges weakly in
$H^{1}$ to some element $\tilde{x}.$ By the compactness of the embedding of
$H^{1}([0,T^{+}],\mathbb{R}^{2N})$ into $C([0,T^{+}],\mathbb{R}^{2N}),$
$\tilde{x}_{k}\rightarrow\tilde{x}$ uniformly on $[0,T^{+}]$. In particular,
$\max_{t\in\lbrack0,T^{+}]}\Vert\tilde{x}(t)\Vert=\lim_{k\rightarrow\infty
}\max_{t\in\lbrack0,T^{+}]}\Vert\tilde{x}_{k}(t)\Vert.$ But
\[
\max_{t\in\lbrack0,T^{+}]}\Vert\tilde{x}_{k}(t)\Vert\geq\Vert\tilde{x}%
_{k}(t_{k}-l_{k}T^{+})\Vert=\left\|  x_{n_{k}}(t_{k})\right\|  \geq
\varepsilon\text{ for all }k,
\]
so $\max_{t\in\lbrack0,T^{+}]}\Vert\tilde{x}(t)\Vert\geq\varepsilon.$ Thus we
see that (2) holds when $t_{k}\rightarrow\infty.$ as $k\rightarrow\infty
.$\newline A similar argument shows that (3) holds if $t_{k}\rightarrow
-\infty.$ as $k\rightarrow\infty,$ completing the proof.
\end{proof}

\begin{theorem}
\label{th4.9} Under the hypotheses (H1) to (H4) and (H$^{\infty})$, suppose
that\newline (1) there is an element $(\lambda,x^{0})\in\mathbb{R}\times
H^{1}$ such that $D_{x}(\lambda,x^{0})\in\Phi_{0}(H^{1},L^{2}),$
and\newline (2)$\{x\in H^{1}:F^{+}(\lambda,x)=0\}=\{x\in H^{1}:F^{-}%
(\lambda,x)=0\}=\{0\}.$\newline Then the restriction of $F(\lambda
,\cdot):H^{1}\rightarrow L^{2}$ to the closed bounded subsets of $H^{1}$ is proper.
\end{theorem}

\begin{proof}
According to Theorem \ref{th4.7} and Lemma
\ref{lem4.6}, it suffices to show that any bounded sequence $\{x_{n}\}$ in
$H^{1} $ such that $\left\|  F(\lambda,x_{n})-y\right\|  _{2}\rightarrow0$ for
some element $y\in L^{2}$ has a weakly convergent subsequence which vanishes
uniformly at infinity. By the boundedness of $\{x_{n}\},$ we may assume
henceforth that $x_{n}\rightharpoonup x$ weakly in $H^{1}$ for some $x\in
H^{1}.$ Furthermore, $\{x_{n}\}$ has at least one of the properties stated in
Lemma \ref{lem4.8} where $T^{\pm}$ are chosen to be periods of $g^{\pm}%
(\cdot,\xi,\lambda)$ as in (H$^{\infty}$)(3).\newline Let us suppose that
$\{x_{n}\}$ has the property (2) of Lemma \ref{lem4.8}. That is to say,
$\widetilde{x_{k}}\rightharpoonup\widetilde{x}$ weakly in $H^{1}$ where
$\widetilde{x_{k}}(t)=x_{n_{k}}(t+l_{k}T^{+})$ and $\widetilde{x}\neq0.$ The
invariance by translation of the Lebesgue measure implies that $\left\|
F(\lambda,x_{n})-y\right\|  _{2}=\left\|  \tau_{l_{k}T^{+}}(F(\lambda
,x_{n})-y)\right\|  _{2}$ so
\[
\left\|  \tau_{l_{k}T^{+}}(F(\lambda,x_{n_{k}}))-\widetilde{y_{k}}\right\|
_{2}\rightarrow0\text{ where }\widetilde{y_{k}}(t)=y(t+l_{k}T^{+})
\]
For any $\omega\in(0,\infty),$ Lemma \ref{lem4.3} shows that
\[
\left\|  \tau_{l_{k}T^{+}}(F(\lambda,x_{n_{k}}))-\tau_{l_{k}T^{+}}%
(F^{+}(\lambda,x_{n_{k}}))\right\|  _{L^{2}(-\omega,\omega)}\rightarrow0\text{
as }k\rightarrow\infty.
\]
Hence
\[
\left\|  \tau_{l_{k}T^{+}}(F^{+}(\lambda,x_{n_{k}}))-\widetilde{y_{k}%
}\right\|  _{L^{2}(-\omega,\omega)}\rightarrow0\text{ as }k\rightarrow\infty.
\]
But
\begin{align*}
\tau_{l_{k}T^{+}}(F^{+}(\lambda,x_{n_{k}}))(t)  & =Jx_{n_{k}}^{\prime}%
(t+l_{k}T^{+})-g^{+}(t+l_{k}T^{+},x_{n_{k}}(t+l_{k}T^{+}),\lambda)\\
& =J\widetilde{x_{k}^{\prime}}(t)-g^{+}(t,\widetilde{x_{k}}(t),\lambda
)=F^{+}(\lambda,\widetilde{x_{k}})(t)
\end{align*}
by the periodicity of $g^{+}.$ Consequently,
\[
\left\|  F^{+}(\lambda,\widetilde{x_{k}})-\widetilde{y_{k}}\right\|
_{L^{2}(-\omega,\omega)}\rightarrow0\text{ as }k\rightarrow\infty,
\]
for all $\omega\in(0,\infty).$ Since the sequence $\{F^{+}(\lambda
,\widetilde{x_{k}})-\widetilde{y_{k}}\}$ is bounded in $L^{2},$ this implies
that $F^{+}(\lambda,\widetilde{x_{k}})-\widetilde{y_{k}}\rightharpoonup0$
weakly in $L^{2}.$ But $\widetilde{y_{k}}\rightharpoonup0$ weakly in $L^{2}$
by Lemma \ref{lem4.4}, so we now have that $F^{+}(\lambda,\widetilde{x_{k}%
})\rightharpoonup0$ weakly in $L^{2}.$ However, the weak sequential continuity
of $F^{+}(\lambda,\cdot):H^{1}\rightarrow L^{2}$ implies that
\[
F^{+}(\lambda,\widetilde{x_{k}})\rightharpoonup F^{+}(\lambda,\widetilde
{x})\text{ weakly in }L^{2},
\]
so we must have $F^{+}(\lambda,\widetilde{x})=0,$ contradicting the hypothesis
(2) of the theorem. This shows that the sequence $\{x_{n}\}$ cannot have the
property (2) of Lemma \ref{lem4.8}.\newline A similar argument excludes the
property (3), completing the proof of the theorem.
\end{proof}

\begin{corollary}
\label{cor4.10} Suppose that (H1) to (H4) and (H$^{\infty})$ are satisfied. An
open interval $\Lambda$ is admissible for $F:\mathbb{R}\times H^{1}\rightarrow
L^{2}$ provided that, for all $\lambda\in\Lambda,$

(1) $D_{x}F(\lambda,0)\in\Phi_{0}(H^{1},L^{2})$ and

(2) $\{x\in H^{1}:F^{+}(\lambda,x)=0\}=\{x\in H^{1}:F^{-}(\lambda,x)=0\}=\{0\}$
\end{corollary}

\begin{proof}
From hypothesis (1) and part (4) of Theorem
\ref{th3.4}, it follows that $D_{x}F(\lambda,x)\in\Phi_{0}(H^{1},L^{2}%
)$\textbf{\ }for all $(\lambda,x)\in\Lambda\times H^{1}.$\newline Let $K$ be a
compact subset of $L^{2}$ and let $W$ be a closed bounded subset of
$\mathbb{R}\times H^{1}$ such that
\[
\inf\Lambda<\inf PW\leq\sup PW<\sup\Lambda.
\]
To show that $F^{-1}(K)\cap W$ is a compact subset of $\mathbb{R}\times
H^{1},$ we consider a sequence $\{(\lambda_{n},x_{n})\}\subset F^{-1}(K)\cap
W.$ Passing to a subsequence, we can suppose that there exist $\lambda
\in\lbrack\inf PW,\sup PW]\subset\Lambda$ and $y\in K$ such that
\[
\lambda_{n}\rightarrow\lambda\text{ and }\left\|  F(\lambda_{n},x_{n}%
)-y\right\|  _{2}\rightarrow0.
\]
But, by part (3) of Theorem \ref{th3.3}, the family of functions
$\{F(\cdot,x_{n})\}_{n\in\mathbb{N}}$ is equicontinuous at $\lambda$ since the
sequence $\{x_{n}\}$ is bounded in $H^{1}.$ It follows from this that
$\left\|  F(\lambda,x_{n})-y\right\|  _{2}\rightarrow0.$ By Theorem
\ref{th4.9} we know that $F(\lambda,\cdot):H^{1}\rightarrow L^{2}$ is proper
on the closed bounded subsets of $H^{1}$ and so there is a subsequence
$\{x_{n_{k}}\}$ of $\{x_{n}\}$ and an element $x\in H^{1}$ such that $\left\|
x_{n_{k}}-x\right\|  \rightarrow0.$ Thus $\{(\lambda_{n_{k}},x_{n_{k}})\}$
converges to $(\lambda,x)$ in $\mathbb{R}\times H^{1}.$ Furthermore
$(\lambda,x)\in W$ since $W$ is closed. This proves that $F^{-1}(K)\cap W$ is
a compact subset of $\mathbb{R}\times H^{1},$ completing the proof that
$\Lambda$ is an admissible interval.
\end{proof}

\begin{theorem}
\label{th4.11} Consider the system (\ref{1.1})(\ref{1.2}) under the
assumptions (H1) to (H4) and (H$^{\infty}).$ Let $\Lambda$ be an open interval
having the following properties.

\begin{itemize}
\item [(1)]For all $\lambda\in\Lambda,D_{x}F(\lambda,0)\in\Phi_{0}(H^{1}%
,L^{2})$.

\item[(2)] For all $\lambda\in\Lambda$,
\[
\{x\in H^{1}:F^{+}(\lambda,x)=0\}=\{x\in H^{1}:F^{-}(\lambda,x)=0\}=\{0\}.
\]

\item[(3)] There is a point $\lambda_{0}\in\Lambda$ such that

\begin{itemize}
\item [(i)]$\dim\ker D_{x}F(\lambda_{0},0)$ is odd,

\item[(ii)] for every $u\in\ker D_{x}F(\lambda_{0},0)\backslash\{0\}$ there is
an element $v\in\ker D_{x}F(\lambda_{0},0)$ such that
\[
\int_{-\infty}^{\infty}\left\langle D_{\lambda}D_{\xi}^{2}H(t,0,\lambda
_{0})u(t),v(t)\right\rangle dt\neq0
\]

\item[(iii)] $\dim\{D_{\lambda}D_{\xi}^{2}H(\cdot,0,\lambda_{0})u:u\in\ker
D_{x}F(\lambda_{0},0)\}=\dim D_{x}F(\lambda_{0},0).$
\end{itemize}
\end{itemize}

\noindent Then a global branch of homoclinic solutions of (\ref{1.1}%
)(\ref{1.2}) bifurcates at $\lambda_{0}$ in the sense of Theorem \ref{th2.3}
with $X=H^{1}$ and $Y=L^{2}.$

\end{theorem}

\begin{theorem}
For some point $t_{0}\in\mathbb{R}$, $\det D_{\lambda}D_{\xi}^{2}%
H(t_{0},0,\lambda_{0})\neq0$.
\end{theorem}

\begin{proof}
In view of Corollary \ref{cor4.10}, we only need to
check that the assumption (3) ensures that the condition (\ref{2.1}) is
satisfied with $A(\lambda)=D_{x}F(\lambda,0).$ By (ii), we have that
$D_{\lambda}D_{\xi}^{2}H(,0,\lambda_{0})u\notin\lbrack\ker D_{x}F(\lambda
_{0},0)]^{\perp}$ for $u\in\ker D_{x}F(\lambda_{0},0)\backslash\{0\}$ and so
$A^{\prime}(\lambda_{0})[\ker D_{x}F(\lambda_{0},0)]\cap$ $\rge A(\lambda
_{0})=\{0\}$ by Theorem \ref{th3.4}(2). Since
\[
A^{\prime}(\lambda_{0})[\ker D_{x}F(\lambda_{0},0)]=\{D_{\lambda}D_{\xi}%
^{2}H(\cdot,0,\lambda_{0})u:u\in\ker D_{x}F(\lambda_{0},0)\}
\]
and codim $\rge A(\lambda_{0})=\dim$ $\ker D_{x}F(\lambda_{0},0)$ by (1), it
follows from (3)(iii) that
\[
A^{\prime}(\lambda_{0})[\ker D_{x}F(\lambda_{0},0)]\oplus\ker A(\lambda
_{0})=L^{2}.
\]
It now follows from (i) and Proposition \ref{prop2.1} that $\pi(D_{x}%
F(\lambda_{0},0)\mid\lbrack\lambda_{0}-\varepsilon,\lambda_{0}+\varepsilon
])=-1$ for some small $\varepsilon>0.$

\noindent\textbf{Remark} \ The condition (3)(iii) is satisfied whenever
\[
(iii)^{\prime}\text{ \ for some point }t_{0}\in\mathbb{R},\det D_{\lambda
}D_{\xi}^{2}H(t_{0},0,\lambda_{0})\neq0.
\]
In fact, if$\{u_{1},...,u_{k}\}$ is a basis for $\ker D_{x}F(\lambda_{0},0)$,
then for every $t\in\mathbb{R},$ the vectors $u_{1}(t),...,u_{k}(t)$ are
linearly independent in $\mathbb{R}^{2N}$ since the functions $u_{1}%
,...,u_{k}$ all satisfy the linear system $Ju^{\prime}(t)=D_{\xi}%
^{2}H(t,0,\lambda_{0})u(t).$ It follows from $(iii)\prime$ that the vectors
\[
D_{\lambda}D_{\xi}^{2}H(t_{0},0,\lambda_{0})u_{1}(t_{0}),\ldots,D_{\lambda
}D_{\xi}^{2}H(t_{0},0,\lambda_{0})u_{k}(t_{0})
\]
are linearly independent in $\mathbb{R}^{2N}$ and hence that $\dim A^{\prime
}(\lambda_{0})[\ker D_{x}F(\lambda_{0},0)]=k$.
\end{proof}

The rest of this paper is devoted to the formulation of explicit conditions on
the Hamiltonian $H$ which enable us to verify the properties (1) to (3) in the
above result.

\section{More explicit criteria}

The first objective is to formulate conditions on the Hamiltonian which ensure
that $D_{x}F(\lambda,0)\in\Phi_{0}(H^{1},L^{2}).$ Recalling that
\[
D_{x}F(\lambda,0)u=Ju^{\prime}-D_{\xi}^{2}H(\cdot,0,\lambda)u\text{ for all
}u\in H^{1},
\]
we set
\[
A_{\lambda}(t)=D_{\xi}^{2}H(t,0,\lambda)\text{ for }t\in\mathbb{R}.
\]
Then, assuming that $(H^{\infty})$ is satisfied, we set
\[
A_{\lambda}^{+}(t)=D_{\xi}g^{+}(t,0,\lambda)\text{ and }A_{\lambda}%
^{-}(t)=D_{\xi}g^{-}(t,0,\lambda).
\]
We observe that $A_{\lambda}(t),A_{\lambda}^{+}(t)$ and $A_{\lambda}^{-}(t)$
are all real, symmetric $2N\times2N-$matrices and that $A_{\lambda}^{+}(t)$
and $A_{\lambda}^{-}(t)$ are periodic in $t.$ Furthermore,
\[
\lim_{t\rightarrow\infty}\{A_{\lambda}(t)-A_{\lambda}^{+}(t)\}=\lim
_{t\rightarrow-\infty}\{A_{\lambda}(t)-A_{\lambda}^{-}(t)\}=0.
\]

\begin{theorem}
\label{th5.1} Suppose that (H1) to (H4) \ and (H$^{\infty})$ are satisfied and
that the periodic, linear Hamiltonian systems
\[
Jx^{\prime}-A_{\lambda}^{+}(t)x=0\text{ and }Jx^{\prime}-A_{\lambda}^{-}(t)x=0
\]
have no characteristic multipliers on the unit circle.\newline Then
$D_{x}F(\lambda,0)\in\Phi_{0}(H^{1},L^{2}).$

\noindent Furthermore,
\begin{equation*}
\ker D_{x}F(\lambda,0)=N(\lambda)\text{ where}
\end{equation*}
\begin{equation*}
N(\lambda)=\{u\in C^{1}(\mathbb{R},\mathbb{R}^{2N}):Ju^{\prime}(t)-A_{\lambda
}(t)u(t)\equiv0\text{ and }\lim_{\left|  t\right|  \rightarrow\infty}u(t)=0\}.
\end{equation*}
\end{theorem}

\noindent\textbf{Remarks }(1) In fact, our proof shows that functions in $\ker
D_{x}F(\lambda,0)$ decay exponentially to zero as $\left|  t\right|
\rightarrow\infty$.

(2) Characteristic multipliers on the unit circle correspond to characteristic
exponents with real part equal to zero.

First we establish the following useful result.

\begin{proposition}
\label{prop5.2} Let $M(t)$ be a real symmetric $2N\times2N-$matrix which
depends continuously and periodically on $t\in\mathbb{R}.$ Suppose that the
linear system
\[
Jx^{\prime}-M(t)x=0
\]
has no characteristic multipliers on the unit circle. Then the linear operator
$L:H^{1}\rightarrow L^{2},$ defined by
\[
Lu=Ju^{\prime}-M(t)u\text{ for all }u\in H^{1},
\]
is an isomorphism.
\end{proposition}

\noindent\textbf{Remark} \ In the case where $M(t)=M$ is constant, there are
no characteristic multipliers on the unit circle precisely when the matrix
$JM$ has no eigenvalues on the imaginary axis.

\begin{proof}
Suppose that $M(t+T)=M(t)$ for all $t\in\mathbb{R}.
$ Recall that a $2N\times2N-$matrix $K$ is symplectic when $K^{T}JK=J$ and
that such matrices are always invertible. By Floquet theory (see for example
Theorem IV-5-11 of \cite{hs}) there exist

\begin{itemize}
\item [(a)]a real symmetric $2N\times2N$-matrix $C$ and

\item[(b)] a real symplectic $2N\times2N$-matrix $P(t)$ for each
$t\in\mathbb{R}$ such that $P_{ij}\in C^{1}(\mathbb{R})$ and $P(t+2T)=P(t)$
for all $t\in\mathbb{R}$, $Ju^{\prime}(t)-M(t)u(t)=P(t)\{Jz^{\prime
}(t)-Cz(t)\}$ for all $t\in\mathbb{R}$ where $u(t)=P(t)z(t)$.
\end{itemize}

Furthermore, the characteristic multipliers of the system
\[
Ju^{\prime
}(t)-M(t)u(t)=0
\]
are the complex numbers $\rho_{1},...,\rho_{2N}$ where
$\rho_{k}=e^{2T\lambda_{k}}$ where $\lambda_{1},...,\lambda_{2N}$ are the
eigenvalues of the matrix $JC.$ For $z:\mathbb{R}\rightarrow\mathbb{R}^{2N},$
let $Wz(t)=P(t)z(t).$ It follows easily from the properties of $P$ that
\begin{align*}
Wz &  \in L^{2}\Longleftrightarrow z\in L^{2}\text{ and also }Wz\in
H^{1}\Longleftrightarrow z\in H^{1},\text{ in fact,}\\
W &  :L^{2}\rightarrow L^{2}\text{ and }W:H^{1}\rightarrow H^{1}\text{ are
isomorphisms.}%
\end{align*}
Thus we can define a bounded linear operator $S:H^{1}\rightarrow L^{2}$ by
setting
\[
Sz=W^{-1}LWz\text{ where }Lu=Ju^{\prime}-M(t)u\text{ for }u\in H^{1}.
\]
Clearly $S:H^{1}\rightarrow L^{2}$ is an isomorphism $\Longleftrightarrow
L:H^{1}\rightarrow L^{2}$ is an isomorphism. But, setting $u=Wz,$%
\[
W^{-1}LWz(t)=P(t)^{-1}\{Ju^{\prime}-M(t)u\}=Jz^{\prime}(t)-Cz(t).
\]
By Corollary 10.2 of \cite{st}, the bounded linear operator $Sz=Jz^{\prime
}(t)-Cz(t)$ is an isomorphism if and only if $\sigma(JC)\cap i\mathbb{R}%
=\emptyset.$ The hypothesis that the system $Ju^{\prime}(t)-M(t)u(t)=0$ has no
characteristic multipliers on the unit circle ensures that indeed
$\sigma(JC)\cap i\mathbb{R}=\emptyset$ and so the proof is complete.
\end{proof}

\noindent\textbf{Proof of Theorem \ref{th5.1}} \ Set $L=D_{x}F(\lambda,0)$ and
recall that
\[
Lx(t)=Jx^{\prime}(t)-A_{\lambda}(t)x(t)\text{ for all }x\in H^{1}.
\]
By Theorem \ref{th3.4}, we know that $L\in B(H^{1},L^{2})$ and it is enough to
show that $\rge L$ is a closed subspace of $L^{2}.$ With this in mind, let
$f\in L^{2}$ and suppose that there exists a sequence $\{f_{n}\}\subset \rge
L$ such that $\left\|  f-f_{n}\right\|  _{2}\rightarrow0.$ Clearly there is a
sequence $\{x_{n}\}\subset H^{1}$ such that $Lx_{n}=f_{n}.$

Let $P$ denote the orthogonal projection of $L^{2}$ onto $\ker L$ and set
$Q=I-P.$ Then $Qx_{n}=x_{n}-Px_{n}\in H^{1}$ and $Qx_{n}\in\lbrack
\ker(L)]^{\perp},$ the orthogonal complement of $\ker L$ in $L^{2}.$ Setting
\[
u_{n}=Qx_{n}\text{ we have that }u_{n}\in H^{1}\cap\lbrack\ker(L)]^{\perp
}\text{ and }Lu_{n}=f_{n}.
\]
Let us prove that the sequence $\{u_{n}\}$ is bounded in $H^{1}.$ For this we
use $S^{\pm}:H^{1}\rightarrow L^{2}$ to denote the bounded linear operators
defined by
\[
S^{\pm}x=Jx^{\prime}-A_{\lambda}^{\pm}(t)x\text{ for all }x\in H^{1}.
\]
By the Proposition \ref{prop5.2} we know that $S^{+}:H^{1}\rightarrow L^{2}$
and $S^{-}:H^{1}\rightarrow L^{2}$ are both isomorphisms and so there exists a
constant $k$ such that
\begin{equation}
\left\|  S^{\pm}x\right\|  _{2}\geq k\left\|  x\right\|  \text{ for all }x\in
H^{1}.\label{5.1}%
\end{equation}
Supposing that $\left\|  u_{n}\right\|  \rightarrow\infty,$ we set
$w_{n}=\frac{u_{n}}{\left\|  u_{n}\right\|  }.$ Then $\{w_{n}\}\subset H^{1}$
with $\left\|  w_{n}\right\|  =1$ for all $n\in\mathbb{N}.$ Passing to a
subsequence, we can suppose that $w_{n}\rightharpoonup w$ weakly in $H^{1},$
and hence that $Lw_{n}\rightharpoonup Lw$ weakly in $L^{2}.$ Furthermore,
\[
Lw_{n}=\frac{Lu_{n}}{\left\|  u_{n}\right\|  }=\frac{f_{n}}{\left\|
u_{n}\right\|  }\text{ and so }\left\|  Lw_{n}\right\|  _{2}=\frac{\left\|
f_{n}\right\|  _{2}}{\left\|  u_{n}\right\|  }\rightarrow0
\]
since $\left\|  f_{n}\right\|  _{2}\rightarrow\left\|  f\right\|  _{2}$ and
$\left\|  u_{n}\right\|  \rightarrow\infty.$ Thus $Lw=0.$ But $w_{n}\in
H^{1}\cap\lbrack\ker(L)]^{\perp}$ for all $n\in\mathbb{N},$ from which it
follows that $w\in H^{1}\cap\lbrack\ker(L)]^{\perp}.$ Thus $w=0$ and
$w_{n}\rightharpoonup0$ weakly in $H^{1}.$ Consequently,
\[
w_{n}\rightarrow0\text{ uniformly on }[-R,R]\text{ for any }R\in(0,\infty).
\]
But, for all $t\in\mathbb{R},$%
\[
Jw_{n}^{\prime}(t)=A_{\lambda}(t)w_{n}(t)+\frac{f_{n}(t)}{\left\|
u_{n}\right\|  }%
\]
and so
\[
\left\|  w_{n}^{\prime}\right\|  _{L^{2}(-R,R)}\leq\sup_{t\in\mathbb{R}%
}\left\|  A_{\lambda}(t)\right\|  \left\|  w_{n}\right\|  _{L^{2}(-R,R)}%
+\frac{\left\|  f_{n}\right\|  _{2}}{\left\|  u_{n}\right\|  },
\]
showing that $\left\|  w_{n}^{\prime}\right\|  _{L^{2}(-R,R)}\rightarrow0$ as
$n\rightarrow\infty,$ for all $R\in(0,\infty).$ In particular,
\[
\left\|  w_{n}\right\|  _{H^{1}(-R,R)}\rightarrow0\text{ as }n\rightarrow
\infty.
\]

Now choose any $\varepsilon>0.$ By (H$^{\infty}),$ there is a constant
$r\in(0,\infty)$ such that
\[
\left|  A_{\lambda}(t)-A_{\lambda}^{+}(t)\right|  \leq\varepsilon\text{ for
all }t\geq r\text{ and }\left|  A_{\lambda}(t)-A_{\lambda}^{-}(t)\right|
\leq\varepsilon\text{ for all }t\leq-r.
\]
There exist a constant $R>r+\frac{1}{\varepsilon}$ and a function $\varphi\in
C^{1}(\mathbb{R})$ such that

\medskip

\noindent$0 \leq\varphi(t)\leq1$ for all $t\in\mathbb{R}$, $\varphi(t)=0$ for
$t\leq r$, $\varphi(t)=1$ for $t\geq R$ and $\left|  \varphi^{\prime
}(t)\right| \leq\varepsilon$ for all $t\in\mathbb{R}$.

\medskip

Consider now the function $z_{n}(t)=\varphi(t)w_{n}(t).$ Clearly $z_{n}\in
H^{1}$ and
\begin{align*}
S^{+}z_{n}(t) &  =\varphi^{\prime}(t)Jw_{n}(t)+\varphi(t)Jw_{n}^{\prime
}(t)-A_{\lambda}^{+}(t)z_{n}(t)\\
&  =\varphi^{\prime}(t)Jw_{n}(t)+\varphi(t)Lw_{n}(t)+\varphi(t)\{A_{\lambda
}(t)-A_{\lambda}^{+}(t)\}w_{n}(t)\\
&  =\varphi^{\prime}(t)Jw_{n}(t)+\varphi(t)\frac{f_{n}(t)}{\left\|
u_{n}\right\|  }+\varphi(t)\{A_{\lambda}(t)-A_{\lambda}^{+}(t)\}w_{n}(t).
\end{align*}
Thus,
\begin{align*}
\left\|  S^{+}z_{n}\right\|  _{2} &  \leq\varepsilon\left\|  w_{n}\right\|
_{2}+\frac{\left\|  f_{n}\right\|  _{2}}{\left\|  u_{n}\right\|  }+\sup_{t\geq
r}\left|  A_{\lambda}(t)-A_{\lambda}^{+}(t)\right|  \left\|  w_{n}\right\|
_{2}\\
&  \leq\varepsilon+\frac{\left\|  f_{n}\right\|  _{2}}{\left\|  u_{n}\right\|
}+\varepsilon
\end{align*}
since $\left\|  w_{n}\right\|  _{2}\leq1.$ Hence, by (\ref{5.1})
\[
\left\|  w_{n}\right\|  _{H^{1}(R,\infty)}=\left\|  z_{n}\right\|
_{H^{1}(R,\infty)}\leq\left\|  z_{n}\right\|  \leq\frac{1}{k}\left( 2\varepsilon
+\frac{\left\|  f_{n}\right\|  _{2}}{\left\|  u_{n}\right\|  }\right).
\]
A similar argument, using $\varphi(-t)w_{n}(t)$ instead of $z_{n},$ shows
that
\[
\left\|  w_{n}\right\|  _{H^{1}(-\infty,-R)}\leq\frac{1}{k} \left( 2\varepsilon
+\frac{\left\|  f_{n}\right\|  _{2}}{\left\|  u_{n}\right\|  }\right).
\]
Finally, we have shown that, for all $n\in\mathbb{N},$%
\begin{align*}
\left\|  w_{n}\right\|  ^{2} &  =\left\|  w_{n}\right\|  _{H^{1}(-\infty
,-R)}^{2}+\left\|  w_{n}\right\|  _{H^{1}(-R,R)}^{2}+\left\|  w_{n}\right\|
_{H^{1}(R,\infty)}^{2}\\
&  \leq\frac{2}{k^{2}}\left( 2\varepsilon+\frac{\left\|  f_{n}\right\|  _{2}%
}{\left\|  u_{n}\right\|  }\right)^{2}+\left\|  w_{n}\right\|  _{H^{1}(-R,R)}^{2}%
\end{align*}
and, letting $n\rightarrow\infty,$%
\[
\limsup_{n\rightarrow\infty}\left\|  w_{n}\right\|  ^{2}\leq\frac{2}{k^{2}%
}\left( 2\varepsilon\right)^{2}%
\]
since $\left\|  w_{n}\right\|  _{H^{1}(-R,R)}^{2}\rightarrow0,\left\|
f_{n}\right\|  _{2}\rightarrow\left\|  f\right\|  _{2}$ and $\left\|
u_{n}\right\|  \rightarrow\infty.$ But $\left\|  w_{n}\right\|  \equiv1$ and
$\varepsilon>0$ can be chosen so that $\frac{2}{k^{2}}\{2\varepsilon\}^{2}<1.$
This contradiction establishes the boundedness of the sequence $\{u_{n}\}$ in
$H^{1}.$

By passing to a subsequence, we can now suppose that $u_{n}\rightharpoonup u$
weakly in $H^{1}$, and consequently that $Lu_{n}\rightharpoonup Lu$ weakly in
$L^{2}.$ However, $Lu_{n}=f_{n}$ and $\left\|  f_{n}-f\right\|  _{2}%
\rightarrow0,$ showing that $Lu=f$. This proves that $\rge L$ is a closed
subspace of $L^{2}$ and we have shown that $D_{x}F(\lambda,0)\in\Phi_{0}%
(H^{1},L^{2}).$

Clearly $\ker D_{x}F(\lambda,0)\subset N(\lambda).$ To establish the equality
we suppose that $u\in N(\lambda)$ and we must show that $u\in H^{1}.$ In fact,
we shall prove that $\left\|  u(t)\right\|  $ decays to zero exponentially as
$\left|  t\right|  \rightarrow\infty.$ This implies that $u\in L^{2}$ and
then, since $A_{\lambda}(t)$ is bounded on $\mathbb{R},$ it follows
immediately that $u^{\prime}\in L^{2}.$ Let us consider the behaviour of
$u(t)$ as $t\rightarrow\infty,$ the case $t\rightarrow-\infty$ being similar.

Using the notation introduced in the proof of Proposition \ref{prop5.2} with
$M=A_{\lambda}^{+},$ we set $u(t)=P(t)z(t)$ and find that $z\in C^{1}%
(\mathbb{R})$ and
\[
P(t)\{Jz^{\prime}(t)-Cz(t)\}=Ju^{\prime}(t)-A_{\lambda}^{+}%
(t)u(t)=\{A_{\lambda}(t)-A_{\lambda}^{+}(t)\}u(t)
\]
where $C$ is a real symmetric matrix and $JC$ has no eigenvalues with zero
real part. Thus,
\[
z^{\prime}(t)=\{-JC+R(t)\}z(t)\text{ where }R(t)=-JP(t)^{-1}\{A_{\lambda
}(t)-A_{\lambda}^{+}(t)\}P(t).
\]
The matrix $R(t)$ depends continuously on $t$ and $\left\|  R(t)\right\|
\rightarrow0$ as $t\rightarrow\infty.$ According to Corollary VII-3-7 of
\cite{hs}, this implies that there is a number $\mu$ such that $\lim
_{t\rightarrow\infty}t^{-1}\log\left\|  z(t)\right\|=\mu$ and $\mu$ is
the real part of an eigenvalue of the matrix $-JC.$ But, since $\lim_{\left|
t\right|  \rightarrow\infty}u(t)=0,$ we have that $\lim_{\left|  t\right|
\rightarrow\infty}z(t)=0$ and consequently $\mu\leq0.$ However, $-JC$ has no
eigenvalues with zero real part and so $\mu<0.$ Thus, for any $\gamma
<-\mu,\lim_{t\rightarrow\infty}e^{\gamma t}\left\|  z(t)\right\|  =0.$ This
establishes the exponential decay of $\left|  u(t)\right|  $ as $t\rightarrow
\infty$ and the proof is complete.

We now turn to the problem of checking the condition
\[
\{x\in H^{1}:F^{+}(\lambda,x)=0\}=\{x\in H^{1}:F^{-}(\lambda,x)=0\}=\{0\}
\]
in Theorem \ref{th4.11}. This amounts to ensuring that certain types of
Hamiltonian system have no solutions which are homoclinic to $0.$

\begin{theorem}
\label{th5.3} Suppose the (H1) to (H4) and (H$^{\infty})$ are satisfied.

(a) If there is a real, symmetric $2N\times2N-$matrix $C$ such that
\[
\left\langle g^{+}(t,\xi,\lambda),JC\xi\right\rangle >0\text{ for all }\xi
\in\mathbb{R}^{2N}\backslash\{0\},
\]
then $\{x\in H^{1}:F^{+}(\lambda,x)=0\}=\{0\}.$

(b) If
\[
g^{+}(t,\xi,\lambda)=A_{\lambda}^{+}(t)\xi\text{ for all }(t,\xi)\in
\mathbb{R}\times\mathbb{R}^{2N},
\]
(that is, $g^{+}(t,\cdot,\lambda)$ is linear), then $\{x\in H^{1}%
:F^{+}(\lambda,x)=0\}=\{0\}.$

(c) If $g^{+}(t,\xi,\lambda)$ is independent of $t$ and $0$ is an isolated
zero of $H^{+}(\cdot,\lambda)$ where $D_{\xi}H^{+}(\cdot,\lambda)=g^{+}%
(\cdot,\cdot,\lambda),$ then $\{x\in H^{1}:F^{+}(\lambda,x)=0\}=\{0\}.$

The same conclusions hold when $g^{+}$ is replaced by $g^{-}.$
\end{theorem}

\begin{proof}
(a) Suppose that $x\in H^{1}$ and $F^{+}%
(\lambda,x)=0.$ Then $x\in C^{1}(\mathbb{R})$ and, for all $t\in\mathbb{R},$%

\[
\frac{d}{dt}\left\langle Cx(t),x(t)\right\rangle =2\left\langle
Cx(t),x^{\prime}(t)\right\rangle =2\left\langle JCx(t),g^{+}(t,x(t),\lambda
)\right\rangle >0
\]
whenever $x(t)\neq0.$ However $\left\langle Cx(t),x(t)\right\rangle
\rightarrow0$ as $t\rightarrow\infty$ and as $t\rightarrow-\infty$ since $x\in
H^{1}.$ Thus we must have that $\left\langle Cx(t),x(t)\right\rangle \equiv0$
and consequently, $\left\langle JCx(t),g^{+}(t,x(t),\lambda)\right\rangle =0$
for all $t\in\mathbb{R}.$ This implies that $x(t)=0$ for all $t\in\mathbb{R}$
as required.

(b) If $x\in H^{1}$ and $F^{+}(\lambda,x)=0,$ it follows that $S^{+}x=0$ in
the notation which was introduced in the proof of Theorem \ref{th5.1}. From
Floquet theory as in the proof of Proposition \ref{prop5.2} with
$M(t)=A_{\lambda}^{+}(t),$ it follows that
\[
S^{+}z=Jz^{\prime}(t)-Cz(t)\text{ where }x(t)=P(t)z(t)\text{ and }z\in H^{1}.
\]
But, as is pointed out at the beginning of Section 10 of \cite{st},
$\ker(S^{+})=\{0\}$ and so $z\equiv0.$ This proves that $x=0.$

(c) If $x\in H^{1}$ and $F^{+}(\lambda,x)=0,$ it follows that $x\in
C^{1}(\mathbb{R})$ and $x$ satisfies the autonomous Hamiltonian system
$Jx^{\prime}(t)=D_{\xi}H(x(t),\lambda).$ Thus $H(x(t),\lambda)$ is constant,
and (\ref{1.2}) implies that $H(x(t),\lambda)=0$ for all $t\in\mathbb{R}.$
Since $0$ is an isolated zero of $H(\cdot,\lambda)$ and $x(t)\rightarrow0$ as
$\left|  t\right|  \rightarrow0,$ it now follows that $x=0$.
\end{proof}

Combining these results we can formulate criteria for admissible intervals
which can be checked in some examples.

\begin{theorem}
\label{th5.4}Suppose that (H1) to (H4) and (H$^{\infty})$ are satisfied. An
open interval $\Lambda$ is admissible provided that, for all $\lambda
\in\Lambda,$ the following conditions are satisfied.

(1) The periodic, linear Hamiltonian systems
\[
Jx^{\prime}-A_{\lambda}^{+}(t)x=0\text{ and }Jx^{\prime}-A_{\lambda}^{-}(t)x=0
\]
have no characteristic multipliers on the unit circle.

(2) The asymptotic limit $g^{+}$ satisfies one of the conditions (a),(b) or
(c) of Theorem \ref{th5.3}

(3) The asymptotic limit $g^{-}$ satisfies one of the conditions (a),(b) or
(c) of Theorem \ref{th5.3}.
\end{theorem}

Finally we can reformulate Theorem \ref{th4.11} as a global bifurcation
theorem concerning the system (\ref{1.1})(\ref{1.2}) with hypotheses only
involving properties of the Hamiltonian.

\begin{theorem}
\label{th5.5}Suppose that (H1) to (H4) and (H$^{\infty})$ are satisfied. An
open interval $\Lambda$ is admissible provided that, for all $\lambda
\in\Lambda,$ the following conditions are satisfied.

(1) The periodic, linear Hamiltonian systems
\[
Jx^{\prime}-A_{\lambda}^{+}(t)x=0\text{ and }Jx^{\prime}-A_{\lambda}^{-}(t)x=0
\]
have no characteristic multipliers on the unit circle.

(2) The asymptotic limit $g^{+}$ satisfies one of the conditions (a),(b) or
(c) of Theorem \ref{th5.3}.

(3) The asymptotic limit $g^{-}$ satisfies one of the conditions (a),(b) or
(c) of Theorem \ref{th5.3}.

(4) There is a point $\lambda_{0}\in\Lambda$ such that

(i) $k=\dim N(\lambda_{0})$ is odd where $N(\lambda)=\{u\in C^{2}%
(\mathbb{R},\mathbb{R}^{2N}):Ju^{\prime}(t)-A_{\lambda}(t)u(t)\equiv0$ and
$\lim_{\left|  t\right|  \rightarrow\infty}u(t)=0\},$

(ii) for every $u\in N(\lambda_{0})\backslash\{0\}$ there exists $v\in
N(\lambda_{0})$ such that
\[
\int_{-\infty}^{\infty}\left\langle T_{\lambda_{0}}(t)u(t),v(t)\right\rangle
dt\neq0\text{ where }T_{\lambda}(t)=D_{\lambda}D_{\xi}^{2}H(t,0,\lambda)\text{
and}%
\]

(iii) $\dim \{T_{\lambda_{0}}(\cdot)u:u\in N(\lambda_{0})\}=k.$

\noindent Then a global branch of homoclinic solutions of (\ref{1.1}%
)(\ref{1.2}) bifurcates at $\lambda_{0}$ in the sense of Theorem \ref{th2.3}
with $X=H^{1}$ and $Y=L^{2}.$
\end{theorem}

Under the hypotheses (H1) to (H4) and (H$^{\infty}),$ all solutions of the
system (\ref{1.1})(\ref{1.2}) decay to zero exponentially fast as $\left|
t\right|  \rightarrow\infty$ and so the system (\ref{1.1})(\ref{1.2}) is
actually equivalent to the equation $F(\lambda,x)=0$ where $F:\mathbb{R}\times
H^{1}\rightarrow L^{2}.$ The exponential decay can be established by a slight
variant of the proof of the second assertion of Theorem \ref{th5.1}.

\begin{theorem}
\label{th5.6} Suppose that (H1) to (H4) and (H$^{\infty})$ are satisfied. If
$x$ is a solution of (\ref{1.1})(\ref{1.2}), there is a constant $\gamma>0$
such that $\lim_{\left|  t\right|  \rightarrow\infty}e^{\gamma\left|
t\right|  }\left|  x(t)\right|  =0.$
\end{theorem}

\begin{proof}
First we observe that
\[
D_{\xi}H(t,\xi,\lambda)=\int_{0}^{1}\frac{d}{d\tau}D_{\xi}H(t,\tau\xi
,\lambda)d\tau=M(t,\xi,\lambda)\xi
\]
where the matrix $M(t,\xi,\lambda)$ is defined by
\[
M(t,\xi,\lambda)=\int_{0}^{1}D_{\xi}^{2}H(t,\tau\xi,\lambda)d\tau.
\]
Thus $x$ satisfies the linear equation
\[
Jx^{\prime}(t)=K(t)x(t)\text{ where }K(t)=M(t,x(t),\lambda).
\]
But,
\begin{align*}
K(t) &  =\int_{0}^{1}\{D_{\xi}^{2}H(t,\tau x(t),\lambda)-D_{\xi}g^{+}(t,\tau
x(t),\lambda)\}d\tau\\
&  +\int_{0}^{1}\{D_{\xi}g^{+}(t,\tau x(t),\lambda)-A_{\lambda}^{+}%
(t)\}d\tau+A_{\lambda}^{+}(t)
\end{align*}
from which it is easy to see that $\left\|  K(t)-A_{\lambda}^{+}(t)\right\|
\rightarrow0$ as $t\rightarrow\infty.$ Now, using the Floquet change of
variables and Corollary VII-3-7 of \cite{hs} as in the proof of Theorem
\ref{th5.1}, we see that $\left|  x(t)\right|  \rightarrow0$ exponentially
fast as $t\rightarrow\infty.$ The behaviour as $t\rightarrow-\infty$ can be
treated in the same way.
\end{proof}

%\newpage

\section{Examples}

Consider the following Hamiltonian,
\begin{align*}
&  H(t,u,v,\lambda)=\\
&  \frac{1}{2}\{v^{2}+\lambda u^{2}+a(t)u^{2}\}+\frac{A(2+\cos t)\left|
u\right|  ^{\sigma+2}}{(\sigma+2)(1+e^{-t})}+\frac{B(2+\cos\omega t)u^{2}%
v^{2}}{2(1+e^{t})}+r(t)Q(u,v,\lambda)
\end{align*}
where $A,B,\sigma,\omega$ are constants with $A\leq0$ and $\sigma>0$,
\[
a,r\in C(\mathbb{R})\text{ with }a
\begin{array}
[c]{c}%
\geq\\
\not \equiv
\end{array}
0\text{ and }\lim_{\left|  t\right|  \rightarrow\infty}a(t)=\lim_{\left|
t\right|  \rightarrow\infty}r(t)=0
\]
and $Q\in C^{3}(\mathbb{R}^{3})$ with
\[
Q(0,0,\lambda)=\partial_{i}Q(0,0,\lambda)=\partial_{j}\partial_{i}%
Q(0,0,\lambda)=0
\]
for all $i,j=1,2$ and all $\lambda\in\mathbb{R}.$

It is easily seen that $H:\mathbb{R}\times\mathbb{R}^{2}\times\mathbb{R}%
\rightarrow\mathbb{R}$ satisfies the conditions (H1) to (H4) and that the
system $Jx^{\prime}(t)=D_{(u,v)}H(t,x(t),\lambda)$ is
\begin{align}
-v^{\prime}(t) &  =\lambda u(t)+a(t)u(t)+\frac{A(2+\cos t)\left|  u(t)\right|
^{\sigma}u(t)}{1+e^{-t}}+\frac{B(2+\cos\omega t)u(t)v(t)^{2}}{1+e^{t}%
}\nonumber\\
&  +r(t)\partial_{u}Q(u(t),v(t),\lambda)\label{6.1}\\
u^{\prime}(t) &  =v(t)+\frac{B(2+\cos\omega t)u(t)^{2}v(t)}{1+e^{t}%
}+r(t)\partial_{v}Q(u(t),v(t),\lambda)\label{6.2}%
\end{align}
The condition (H$^{\infty}$) is also satisfied with
\begin{align*}
H^{+}(t,u,v,\lambda) &  =\frac{1}{2}\{v^{2}+\lambda u^{2}\}+\frac{A(2+\cos
t)\left|  u\right|  ^{\sigma+2}}{\sigma+2},\\
H^{-}(t,u,v,\lambda) &  =\frac{1}{2}\{v^{2}+\lambda u^{2}\}+\frac
{B(2+\cos\omega t)u^{2}v^{2}}{2},\\
g^{+}(t,u,v,\lambda) &  =(\lambda u+A(2+\cos t)\left|  u\right|  ^{\sigma
}u,v)\\
g^{-}(t,u,v,\lambda) &  =(\lambda u+B(2+\cos\omega t)uv^{2},v+B(2+\cos\omega
t)u^{2}v).
\end{align*}
Note that
\begin{align*}
A_{\lambda}(t) &  =D_{(u,v)}^{2}H(t,0,0,\lambda)=\left[
\begin{array}
[c]{cc}%
\lambda+a(t) & 0\\
0 & 1
\end{array}
\right]  \text{ and}\\
A_{\lambda}^{\pm}(t) &  =D_{(u,v)}g^{\pm}(t,0,0,\lambda)=\left[
\begin{array}
[c]{cc}%
\lambda & 0\\
0 & 1
\end{array}
\right]  .
\end{align*}
Thus $A_{\lambda}^{\pm}(t)$ is independent of $t$ and the spectrum of the
matrix $J\left[
\begin{array}
[c]{cc}%
\lambda & 0\\
0 & 1
\end{array}
\right]  $ is
\[
\{\pm i\sqrt{\lambda}\}\text{ for }\lambda>0,\qquad\{0\}\text{ for }%
\lambda=0\text{ and }\{\pm\sqrt{\left|  \lambda\right|  }\}\text{ for }%
\lambda<0.
\]
It follows from Theorem \ref{th5.1} that $D_{x}F(\lambda,0)\in\Phi_{0}%
(H^{1},L^{2})$ for all $\lambda<0.$

Setting $C=\left[
\begin{array}
[c]{cc}%
0 & 1\\
1 & 0
\end{array}
\right]  ,$ we find that $JC(u,v)^{T}=(-u,v)^{T}$ and hence that
\[
\left\langle g^{+}(t,u,v,\lambda),JC(u,v)^{T}\right\rangle =-\lambda
u(t)^{2}-A(2+\cos t)\left|  u\right|  ^{\sigma+2}+v^{2}>0
\]
for all $(u,v)\in\mathbb{R}^{2}\backslash\{(0,0)\}$ provided that $\lambda<0.$
Similarly,
\[
\left\langle g^{-}(t,u,v,\lambda),JC(u,v)^{T}\right\rangle =-\lambda
u(t)^{2}+v^{2}>0
\]
for all $(u,v)\in\mathbb{R}^{2}\backslash\{(0,0)\}$ provided that $\lambda<0.$

By Theorem \ref{th5.4} we now have that $(-\infty,0)$ is an admissible
interval for the system (\ref{6.1})(\ref{6.2}). In particular, $D_{x}%
F(\lambda,0):H^{1}\rightarrow L^{2}$ is an isomorphism for $\lambda<0, $
except at the values of $\lambda$ for which $\ker D_{x}F(\lambda,0)\neq\{0\}.$
But, $x\in\ker D_{x}F(\lambda,0)\backslash\{0\}$ means that $x$ is a
homoclinic solution of the linear system
\[%
\begin{array}
[c]{cc}%
-x_{2}^{\prime}(t)= & \{\lambda+a(t)\}x_{1}(t)\\
x_{1}^{\prime}(t)= & x_{2}(t)
\end{array}
\]
which is equivalent to the second order equation
\[
x_{1}^{\prime\prime}(t)=-\{\lambda+a(t)\}x_{1}(t).
\]
Under our hypotheses on the coefficient $a,$ there is always at least one
value of $\lambda$ in the interval $(-\infty,0)$ for which this equation has a
homoclinic solution. Setting
\[
\lambda_{0}=\inf\{\textstyle\int_\mathbb{R}\varphi^{\prime}(t)^{2}%
-a(t)\varphi(t)^{2}dt:\varphi\in H^{1}(\mathbb{R})\text{ with 
}\textstyle\int_{\mathbb{R}}\varphi(t)^{2}dt=1\},
\]
it is well-known (see \cite{LL}, theorem 11.5) that $\lambda_{0}\in(-\infty,0)$ and that there 
exists an
element $\varphi_{0}\in H^{1}(\mathbb{R})$ such that
\[
\varphi_{0}(t)>0\text{ for all }t\in\mathbb{R}\text{ and }\lambda_{0}%
\int_{\mathbb{R}}\varphi(t)_{0}^{2}dt=\int_{\mathbb{R}-}\varphi
_{0}^{\prime}(t)^{2}-a(t)\varphi_{0}(t)^{2}dt.
\]
Furthermore, $\varphi_{0}\in H^{2}(\mathbb{R})\cap C^{2}(\mathbb{R})$ and
satisfies the equation
\[
x_{1}^{\prime\prime}(t)=-\{\lambda_{0}+a(t)\}x_{1}(t).
\]
Setting $x_{0}=(\varphi_{0},\varphi_{0}^{\prime}),$ we find that $\ker
D_{x}F(\lambda_{0},0)=\operatorname{span}\{x_{0}\}.$

\noindent Finally we observe that
\[
D_{\lambda}D_{(u,v)}^{2}H(t,0,0,\lambda)=\left[
\begin{array}
[c]{cc}%
1 & 0\\
0 & 0
\end{array}
\right]
\]
and so, in the notation of Theorem \ref{th5.5},
\[
T_{\lambda_{0}}(t)x_{0}(t)=(\varphi_{0},0).
\]
Thus all the hypotheses of Theorem \ref{th5.5} are satisfied by the system
(\ref{6.1})(\ref{6.2}).

\section{Appendix}

\noindent\textbf{Proof of Theorem \ref{th3.3} }(1) It is sufficient to prove
that $D_{x}F$ and $D_{\lambda}F$ exist and are continuous on $\mathbb{R}\times
H^{1}.$\newline For $D_{x}F(\lambda,x),$ we consider $\lambda\in\mathbb{R}$
and $x,y\in H^{1}$ with $(\lambda,x)$ fixed$.$ Then
\begin{align*}
& F(\lambda,x+y)-F(\lambda,x)-\{Jy^{\prime}-M(\lambda,x)y\}\\
& =-D_{\xi}H(\cdot,x+y,\lambda)+D_{\xi}H(\cdot,x,\lambda)+D_{\xi}^{2}%
H(\cdot,x,\lambda)y\\
& =\int_{0}^{1}D_{\xi}^{2}H(\cdot,x,\lambda)y-\frac{d}{ds}D_{\xi}%
H(\cdot,x+sy,\lambda)ds\\
& =\int_{0}^{1}\{D_{\xi}^{2}H(\cdot,x,\lambda)-D_{\xi}^{2}H(\cdot
,x+sy,\lambda)\}yds.
\end{align*}
Hence, for all $t\in\mathbb{R},$%
\begin{align*}
& \left\|  F(\lambda,x+y)(t)-F(\lambda,x)(t)-\{Jy^{\prime}-M(\lambda
,x)y\}(t)\right\| \\
& \leq\left\|  y(t)\right\|  \sup_{t\in\mathbb{R}}\int_{0}^{1}\left\|  D_{\xi
}^{2}H(t,x(t),\lambda)-D_{\xi}^{2}H(t,x(t)+sy(t),\lambda)\right\|  ds
\end{align*}
and so
\begin{align*}
& \left\|  F(\lambda,x+y)-F(\lambda,x)-\{Jy^{\prime}-M(\lambda,x)y\}\right\|
_{2}\\
& \leq\left\|  y\right\|  _{2}\sup_{t\in\mathbb{R}}\int_{0}^{1}\left\|
D_{\xi}^{2}H(t,x(t),\lambda)-D_{\xi}^{2}H(t,x(t)+sy(t),\lambda)\right\|  ds.
\end{align*}
Recalling that $H^{1}$ is continuously embedded in $C_{d},$ we observe that
there is a compact subset $K$ of $\mathbb{R}^{2N}$ such that $x(t)$ and
$y(t)\in K$ for all $t\in\mathbb{R}$ and all $y\in H^{1}$ such that $\left\|
y\right\|  \leq1.$ Since $D_{\xi}^{2}H(\cdot,\cdot,\lambda)$ is a $C_{\xi
}^{0}-$bundle map by (H3), it follows that
\[
\sup_{t\in\mathbb{R}}\int_{0}^{1}\left\|  D_{\xi}^{2}H(t,x(t),\lambda)-D_{\xi
}^{2}H(t,x(t)+sy(t),\lambda)\right\|  ds\rightarrow0\text{ as }\left\|
y\right\|  \rightarrow0.
\]
This proves that $D_{x}F(\lambda,x)y$ exists and is equal to $Jy^{\prime
}-M(\lambda,x)y.$ For the continuity of $D_{x}F,$ we consider $(\lambda
,x),(\mu,z)\in\mathbb{R}\times H^{1}$ and $y\in H^{1}$ with $(\lambda,x)$
fixed. Then
\begin{align*}
\{D_{x}F(\lambda,x)-D_{x}F(\mu,z)\}y  & =\{M(\mu,z)-M(\lambda,x)\}y\\
& =\{D_{\xi}^{2}H(\cdot,z,\mu)-D_{\xi}^{2}H(\cdot,x,\lambda)\}y
\end{align*}
and hence
\[
\left\|  \{ D_{x}F(\lambda,x)-D_{x}F(\mu,z)\}y\right\|  _{2}\leq\left\|
y\right\|  _{2}\sup_{t\in\mathbb{R}}\left\|  D_{\xi}^{2}H(t,z(t),\mu)-D_{\xi
}^{2}H(t,x(t),\lambda)\right\|  .
\]
Thus $\left\|  D_{x}F(\lambda,x)-D_{x}F(\mu,z)\right\|  $ in $B(H^{1},L^{2})$
is bounded above by
\[
\sup_{t\in\mathbb{R}}\left\|  D_{\xi}^{2}H(t,z(t),\mu)-D_{\xi}^{2}%
H(t,x(t),\lambda)\right\|  .
\]
But
\begin{align*}
& D_{\xi}^{2}H(t,z(t),\mu)-D_{\xi}^{2}H(t,x(t),\lambda)\\
& =\int_{0}^{1}\frac{d}{ds}D_{\xi}^{2}H(t,z(t),s\mu+(1-s)\lambda)ds+D_{\xi
}^{2}H(t,z(t),\lambda)-D_{\xi}^{2}H(t,x(t),\lambda)\\
& =\int_{0}^{1}D_{\lambda}D_{\xi}^{2}H(t,z(t),s\mu+(1-s)\lambda)ds(\mu
-\lambda)+D_{\xi}^{2}H(t,z(t),\lambda)-D_{\xi}^{2}H(t,x(t),\lambda)
\end{align*}
There is a constant $K$ such that $\left\|  (z(t),s\mu+(1-s)\lambda)\right\|
\leq K$ for all $(\mu,z)\in\mathbb{R}\times H^{1}$ with $\left\|
(\mu,z)-(\lambda,x)\right\|  \leq1,$ and hence by part (i) of Lemma
\ref{lem3.2}, there is a constant $C(K)$ such that
\[
\left\|  D_{\lambda}D_{\xi}^{2}H(t,z(t),s\mu+(1-s)\lambda)\right\|  \leq
C(K)\text{ for all }t\in\mathbb{R}\text{ and all }s\in\lbrack0,1].
\]
Thus
\[
\sup_{t\in\mathbb{R}}\int_{0}^{1}\left\|  D_{\lambda}D_{\xi}^{2}%
H(t,z(t),s\mu+(1-s)\lambda)\right\|  ds\leq C(K)
\]
and so
\begin{align*}
& \sup_{t\in\mathbb{R}}\left\|  D_{\xi}^{2}H(t,z(t),\mu)-D_{\xi}%
^{2}H(t,x(t),\lambda)\right\| \\
& \leq C(K)\left|  \mu-\lambda\right|  +\sup_{t\in\mathbb{R}}\left\|  D_{\xi
}^{2}H(t,z(t),\lambda)-D_{\xi}^{2}H(t,x(t),\lambda)\right\|  .
\end{align*}
As above,
\[
\sup_{t\in\mathbb{R}}\left\|  D_{\xi}^{2}H(t,z(t),\lambda)-D_{\xi}%
^{2}H(t,x(t),\lambda)\right\|  \rightarrow0\text{ as }\left\|  z-x\right\|
\rightarrow0
\]
since $D_{\xi}^{2}H(\cdot,\cdot,\lambda)$ is a $C_{\xi}^{0}-$bundle map. Thus
\[
\sup_{t\in\mathbb{R}}\left\|  D_{\xi}^{2}H(t,z(t),\mu)-D_{\xi}^{2}%
H(t,x(t),\lambda)\right\|  \rightarrow0
\]
as $(\mu,z)\rightarrow(\lambda,x)$ in $\mathbb{R}\times H^{1},$ establishing
the continuity of

\noindent$D_{x}F:$ $\mathbb{R}\times H^{1}\rightarrow B(H^{1},L^{2})$ at
$(\lambda,x).$

For the differentiability with respect to $\lambda,$ we consider $\lambda
,\tau\in\mathbb{R}$ and $x\in H^{1}.$ Then
\begin{align*}
& F(\lambda+\tau,x)-F(\lambda,x)+\tau D_{\lambda}D_{\xi}H(\cdot,x,\lambda)\\
& =-D_{\xi}H(\cdot,x,\lambda+\tau)+D_{\xi}H(\cdot,x,\lambda)+\tau D_{\lambda
}D_{\xi}H(\cdot,x,\lambda)\\
& =\int_{0}^{1}\{\tau D_{\lambda}D_{\xi}H(\cdot,x,\lambda)-\frac{d}{ds}D_{\xi
}H(\cdot,x,\lambda+s\tau)\}ds\\
& =\tau\int_{0}^{1}\{D_{\lambda}D_{\xi}H(\cdot,x,\lambda)-D_{\lambda}D_{\xi
}H(\cdot,x,\lambda+s\tau)\}ds.
\end{align*}
By (H1), $D_{\lambda}D_{\xi}H(\cdot,0,\lambda)\equiv0$ and so
\begin{align*}
& D_{\lambda}D_{\xi}H(\cdot,x,\lambda)-D_{\lambda}D_{\xi}H(\cdot
,x,\lambda+s\tau)\\
& =\int_{0}^{1}\frac{d}{d\sigma}\{D_{\lambda}D_{\xi}H(\cdot,\sigma
x,\lambda)-D_{\lambda}D_{\xi}H(\cdot,\sigma x,\lambda+s\tau)\}d\sigma\\
& =\int_{0}^{1}\{D_{\lambda}D_{\xi}^{2}H(\cdot,\sigma x,\lambda)-D_{\lambda
}D_{\xi}^{2}H(\cdot,\sigma x,\lambda+s\tau)\}xd\sigma.
\end{align*}
Thus
\begin{align*}
& \left\|  \frac{F(\lambda+\tau,x)-F(\lambda,x)}{\tau}+D_{\lambda}D_{\xi
}H(\cdot,x,\lambda)\right\|  _{2}\\
& \leq\left\|  x\right\|  _{2}\sup_{t\in\mathbb{R}}\int_{0}^{1}\int_{0}%
^{1}\left\|  D_{\lambda}D_{\xi}^{2}H(t,\sigma x(t),\lambda)-D_{\lambda}D_{\xi
}^{2}H(t,\sigma x(t),\lambda+s\tau)\right\|  d\sigma\, ds.
\end{align*}
Recalling that $D_{\lambda}D_{\xi}^{2}H$ is a $C_{(\xi,\lambda)}^{0}-$bundle
map, we easily deduce from this that $D_{\lambda}F(\lambda,x)$ exists and is
equal to $-D_{\lambda}D_{\xi}H(\cdot,x,\lambda)$ for all $(\lambda
,x)\in\mathbb{R}\times H^{1}.$ For the continuity of $D_{\lambda}F,$ we
consider $(\lambda,x),(\mu,z)\in\mathbb{R}\times H^{1}$ with $(\lambda,x)$
fixed. Then, using (H1),
\begin{align*}
D_{\lambda}F(\lambda,x)-D_{\lambda}F(\mu,z)  & =D_{\lambda}D_{\xi}%
H(\cdot,z,\mu)-D_{\lambda}D_{\xi}H(\cdot,x,\lambda)\\
& =\int_{0}^{1}\frac{d}{ds}\{D_{\lambda}D_{\xi}H(\cdot,sz,\mu)-D_{\lambda
}D_{\xi}H(\cdot,sx,\lambda)\}ds\\
& =\int_{0}^{1}\{D_{\lambda}D_{\xi}^{2}H(\cdot,sz,\mu)z-D_{\lambda}D_{\xi}%
^{2}H(\cdot,sx,\lambda)x\}ds
\end{align*}

Hence
%\begin{align*}
\begin{multline*}
 D_{\lambda}F(\lambda,x)-D_{\lambda}F(\mu,z)=\\
=\int_{0}^{1}\{D_{\lambda}D_{\xi}^{2}H(\cdot,sz,\mu)z-D_{\lambda}D_{\xi}%
^{2}H(\cdot,sx,\lambda)z\}ds\\
+\int_{0}^{1}\{D_{\lambda}D_{\xi}^{2}H(\cdot,sx,\lambda)z-D_{\lambda}D_{\xi
}^{2}H(\cdot,sx,\lambda)x\}ds
%\end{align*}
\end{multline*}
and so
%\begin{align*}
\begin{multline*}
\left\|  D_{\lambda}F(\lambda,x)-D_{\lambda}F(\mu,z)\right\|  _{2}\\
\leq\left\|  z\right\|  _{2}\sup_{t\in\mathbb{R}}\int_{0}^{1}\left\|
D_{\lambda}D_{\xi}^{2}H(t,sz(t),\mu)-D_{\lambda}D_{\xi}^{2}H(t,sx(t),\lambda
)\right\|  ds\\
+\left\|  z-x\right\|  _{2}\sup_{t\in\mathbb{R}}\int_{0}^{1}\left\|
D_{\lambda}D_{\xi}^{2}H(t,sx(t),\lambda)\right\|  ds
%\end{align*}
\end{multline*}

From (H3), it follows that
\[
\sup_{t\in\mathbb{R}}\int_{0}^{1}\left\|  D_{\lambda}D_{\xi}^{2}%
H(t,sz(t),\mu)-D_{\lambda}D_{\xi}^{2}H(t,sx(t),\lambda)\right\|
ds\rightarrow0
\]
as $(z,\mu)\rightarrow(x,\lambda)$ in $\mathbb{R}\times H^{1},$ and it follows
from Lemma \ref{lem3.2}(i) that there is a constant $C$ such that
\[
\sup_{t\in\mathbb{R}}\int_{0}^{1}\left\|  D_{\lambda}D_{\xi}^{2}%
H(\cdot,sx,\lambda)\right\|  ds\leq C.
\]
Therefore
\[
\left\|  D_{\lambda}F(\lambda,x)-D_{\lambda}F(\mu,z)\right\|  _{2}%
\rightarrow0\text{ as }(z,\mu)\rightarrow(x,\lambda)\text{ in }\mathbb{R}%
\times H^{1},
\]
proving that $D_{\lambda}F:\mathbb{R}\times H^{1}\rightarrow L^{2}$ is
continuous at $(\lambda,x).$

(2) To prove that $D_{x}F(\lambda,0)$ is differentiable with respect to
$\lambda,$ we consider $\lambda,\tau\in\mathbb{R}$ and $y\in H^{1}.$ Then
\begin{align*}
\{D_{x}F(\lambda+\tau,0)-D_{x}F(\lambda,0)\}y  & =-\{D_{\xi}^{2}%
H(\cdot,0,\lambda+\tau)-D_{\xi}^{2}H(\cdot,0,\lambda)\}y\\
& =-\tau\int_{0}^{1}D_{\lambda}D_{\xi}^{2}H(\cdot,0,\lambda+s\tau)yds
\end{align*}
and hence
%\begin{align*}
\begin{multline*}
\left\{ \frac{D_{x}F(\lambda+\tau,0)-D_{x}F(\lambda,0)}{\tau}+D_{\lambda}D_{\xi
}^{2}H(\cdot,0,\lambda)\right\} y\\
=\int_{0}^{1}\{D_{\lambda}D_{\xi}^{2}H(\cdot,0,\lambda)-D_{\lambda}D_{\xi
}^{2}H(\cdot,0,\lambda+s\tau)\}yds.
%\end{align*}
\end{multline*}
Thus
%\begin{align*}
\begin{multline*}
\left\|  \left\{ \frac{D_{x}F(\lambda+\tau,0)-D_{x}F(\lambda,0)}{\tau}+D_{\lambda
}D_{\xi}^{2}H(\cdot,0,\lambda)\right\} y\right\|  _{2}\\
\leq\left\|  y\right\|  _{2}\sup_{t\in\mathbb{R}}\int_{0}^{1}\left\|
D_{\lambda}D_{\xi}^{2}H(\cdot,0,\lambda)-D_{\lambda}D_{\xi}^{2}H(\cdot
,0,\lambda+s\tau)\right\|  ds
%\end{align*}
\end{multline*}
and, since $H^{1}$ is continuously embedded in $L^{2},$ this shows that
%\begin{align*}
\begin{multline*}
\left\|  \frac{D_{x}F(\lambda+\tau,0)-D_{x}F(\lambda,0)}{\tau}+D_{\lambda
}D_{\xi}^{2}H(\cdot,0,\lambda)[\cdot]\right\|  _{B(H^{1},L^{2})}\\
\leq\sup_{t\in\mathbb{R}}\int_{0}^{1}\left\|  
D_{\lambda}D_{\xi}^{2}H(t,0,\lambda)-D_{\lambda}D_{\xi}^{2}H(t,0,\lambda+s\tau)\right\|  ds
%\end{align*}
\end{multline*}
Using (H3), it follows that $D_{\lambda}D_{x}F(\lambda,0)$ exists and is equal
to multiplication by $-D_{\lambda}D_{\xi}^{2}H(\cdot,0,\lambda).$ The
continuity of $D_{\lambda}D_{x}F(\cdot,0):\mathbb{R}\rightarrow B(H^{1}%
,L^{2})$ again follows from (H3). Indeed,
\begin{align*}
\left\|  \{D_{\lambda}D_{x}F(\lambda,0)-D_{\lambda}D_{x}F(\mu,0)\}y\right\|
_{2}  & =\left\|  \{D_{\lambda}D_{\xi}^{2}H(\cdot,0,\mu)-D_{\lambda}D_{\xi
}^{2}H(\cdot,0,\lambda)\}y\right\|  _{2}\\
& \leq\left\|  y\right\|  _{2}\sup_{t\in\mathbb{R}}\left\|  D_{\lambda}D_{\xi
}^{2}H(t,0,\mu)-D_{\lambda}D_{\xi}^{2}H(t,0,\lambda)\right\|
\end{align*}
so that
\[
\left\|  D_{\lambda}D_{x}F(\lambda,0)-D_{\lambda}D_{x}F(\mu,0)\right\|
_{B(H^{1},L^{2})}\leq\sup_{t\in\mathbb{R}}\left\|  D_{\lambda}D_{\xi}%
^{2}H(t,0,\mu)-D_{\lambda}D_{\xi}^{2}H(t,0,\lambda)\right\|  .
\]

(3) Consider $u\in W$ and $\lambda,\mu\in\mathbb{R}$ with $\lambda$ fixed.
Then
\[
\left\|  F(\lambda,u)-F(\mu,u)\right\|  _{2}=\left\|  D_{\xi}H(\cdot
,u,\lambda)-D_{\xi}H(\cdot,u,\mu)\right\|  _{2}%
\]
and
\begin{align*}
D_{\xi}H(\cdot,u,\lambda)-D_{\xi}H(\cdot,u,\mu)  & =(\lambda-\mu)\int_{0}%
^{1}D_{\lambda}D_{\xi}H(\cdot,u,s\lambda+(1-s)\mu)ds\\
& =(\lambda-\mu)\int_{0}^{1}\int_{0}^{1}D_{\lambda}D_{\xi}^{2}H(\cdot,\sigma
u,s\lambda+(1-s)\mu)udsd\sigma.
\end{align*}
since $D_{\lambda}D_{\xi}H(\cdot,0,s\lambda+(1-s)\mu)\equiv0$ by (H1). Let
$L=\sup\{\left\|  x\right\|  _{\infty}:x\in W\}$. Since $W$ is bounded in
$H^{1}$ which is continuously embedded in $C_{d}$, it follows that $L<\infty.$
Set $K=\sqrt{L+\left|  \lambda\right|  +1}.$ By Lemma \ref{lem3.2}(i) there is
a constant $C(K)$ such that
\[
\left\|  D_{\lambda}D_{\xi}^{2}H(t,\xi,\gamma)\right\|  \leq C(K)\text{ for
all }t\in\mathbb{R}\text{ and }\left\|  (\xi,\gamma)\right\|  \leq K.
\]
It follows that, for all $u\in W$ and $\mu\in\lbrack\lambda-1,\lambda+1],$%
\[
\left\|  F(\lambda,u)-F(\mu,u)\right\|  _{2}\leq\left|  \lambda-\mu\right|
\left\|  u\right\|  _{2}C(K),
\]
proving that the family $\{F(\cdot,u)\}_{u\in W}$ is equicontinuous at
$\lambda.$

(4) We fix $(\lambda,x)\in\mathbb{R}\times H^{1}$ and consider a sequence
$\{x_{n}\}$ which converges weakly to $x$ in $H^{1}.$ In particular,\ $\{x_{n}%
\}$ is bounded in $H^{1},$ and recalling the remarks following (\ref{3.3}), we
have that $\{F(\lambda,x_{n})\}$ is a bounded sequence in $L^{2}.$ Thus it is
enough to prove that
\[
\left\langle F(\lambda,x_{n})-F(\lambda,x_{n}),\varphi\right\rangle
_{2}\rightarrow0\text{ for all }\varphi\in C_{0}^{\infty}(\mathbb{R}%
,\mathbb{R}^{2N})=C_{0}^{\infty}%
\]
where $\left\langle \cdot,\cdot\right\rangle _{2}$ denotes the usual scalar
product in $L^{2}.$\newline Given $\varphi\in C_{0}^{\infty},$ let $R>0$ be
such that $\varphi(t)=0$ for all $t\notin\lbrack-R,R].$ Furthermore, since
$H^{1}$ is continuously embedded in $C_{d},$ there is a constant $K$ such
that
\[
\left\|  x_{n}\right\|  _{\infty}\leq K\text{ for all }n\in\mathbb{N}\text{
and also }\left\|  x\right\|  _{\infty}\leq K.
\]
Since $D_{\xi}H(\cdot,\cdot,\lambda)$ is uniformly continuous on $[-R,R]\times
B(0,2K)\subset\mathbb{R}\times\mathbb{R}^{2N},$ it follows that
\[
\left\langle D_{\xi}H(\cdot,x_{n},\lambda)-D_{\xi}H(\cdot,x,\lambda
),\varphi\right\rangle \rightarrow0\text{ as }n\rightarrow\infty
\]
since $x_{n}\rightarrow x$ uniformly on $[-R,R].$ Furthermore,
\[
\left\langle Jx_{n}^{\prime}-Jx^{\prime},\varphi\right\rangle _{2}%
=\left\langle x_{n}-x,J\varphi^{\prime}\right\rangle _{2}\rightarrow0
\]
again by the uniform convergence of $x_{n}$ to $x$ on $[-R,R].$ This completes
the proof.

\bigskip

\noindent\textbf{Proof of Theorem \ref{th4.1} }By (H$^{\infty})$ we have
that $D_{\xi}g^{+}(\cdot,\cdot,\lambda)$ is a $C_{\xi}^{0}-$bundle map and
\[
g^{+}(t,\xi,\lambda)=\int_{0}^{1}D_{\xi}g^{+}(t,s\xi,\lambda)\xi ds.
\]
As in the proof of Lemma \ref{lem3.2}(ii), for any $K>0,$ there exists a
constant $C(\lambda,K)$ such that
\[
\left\|  D_{\xi}g^{+}(t,\xi,\lambda)\right\|  \leq C(\lambda,K)
\]
for all $t\in\mathbb{R}$ and all $\xi\in\mathbb{R}^{2N}$ such that $\left\|
\xi\right\|  \leq K.$ Thus,
\[
\left\|  g^{+}(t,\xi,\lambda)\right\|  \leq C(\lambda,K)\left\|  \xi\right\|
\]
for all $t\in\mathbb{R}$ and all $\xi\in\mathbb{R}^{2N}$ such that $\left\|
\xi\right\|  \leq K.$ If $W$ is a bounded subset of $H^{1},$ there is a
constant $K$ such that $\left\|  x\right\|  _{\infty}\leq K$ for all $x\in W,$
and consequently,
\[
\left\|  g^{+}(t,x(t),\lambda)\right\|  \leq C(\lambda,K)\left\|
x(t)\right\|
\]
for all $t\in\mathbb{R}$ and $x\in W.$ Thus $g^{+}(\cdot,x,\lambda)\in L^{2}$
and
\[
\left\|  g^{+}(\cdot,x,\lambda)\right\|  _{2}\leq C(\lambda,K)\left\|
x\right\|  _{2}\text{ for all }x\in W.
\]
This proves that $h^{+}(\lambda,\cdot)$ maps $H^{1}$ into $L^{2}$ and is bounded.

The weakly sequential continuity of $F(\lambda,\cdot):H^{1}\rightarrow L^{2}$
follows from this and the continuity of $g^{+}(\cdot,\cdot,\lambda
):\mathbb{R}\times\mathbb{R}^{2N}\rightarrow\mathbb{R}^{2N}$ exactly as in the
proof of Theorem \ref{th3.3}(4).

\bigskip

\noindent
{\bf Acknowledgement. }  This paper was begun when the first author was visiting the EPFL 
in Lausanne. He would like to thank all the staff for providing a warm environment. Moreover, 
he wishes to thank prof. J.~Pejsachowicz for many interesting discussions.

\end{document}